\documentclass[12pt]{amsart}
\usepackage[cp1251]{inputenc}
\usepackage[english,russian]{babel}
\usepackage{amsmath}
\usepackage{amssymb}
\usepackage{amsfonts}
\usepackage{graphicx}
\newtheorem{lemma}{Лемма}
\newtheorem{theorem}{Теорема}

\newtheorem{problem}{Задача}

\begin{document}

 УДК{517.984.50}
 %\udcs \thispagestyle

\title[ Ряды Фурье оператора ротор]
{Ряды Фурье оператора ротор и пространства Соболева}
\author{\bf {Р.С.~Сакс}}
%\address{Ромэн Семенович Сакс Краевые задачи для
%\newline\hphantom{iii} Институт математики c ВЦ УНЦ  РАН,
%\newline\hphantom{iii} ул.Чернышевского, 112,
%\newline\hphantom{iii} 450077, г. Уфа, Россия}
%\email{romen-saks@yandex.ru}

%©   Р.С.Сакс 2012-09-12

\maketitle
 Аннотация. {\small В работе описаны
   свойства операторов ротор и  градиент дивергенции
    в   области $G$ трехмерного пространства. Обсуждается само-сопряженность этих
   операторов в подпространствах $\mathbf{L}_{2}(G)$ и базисность системы собственных функций. Выписаны явные
   формулы для решения краевых задач в шаре и условия разложимости вектор-функций
   в ряды Фурье по собственным функциям ротора и  градиента дивергенции.}
     % , их спектральные разложения и
     %краевые задачи.%автор изучает%  решены явно произвольной%ограниченной
\section {Введение и основные результаты}

\subsection {Основные пространства} В статье мы рассматриваем линейные
пространства над полем $\mathbb{C}$ комплексных чисел. Через
$\mathbf{L}_{2}(G)$ обозначаем пространство Лебега вектор-функций,
квадратично интегрируемых в $G$    с внутренним произведением
$(\mathbf {u},\mathbf {v})= \int_G \,\mathbf
{u}\cdot\overline{\mathbf {v}}\,d\,\mathbf {x}$  и нормой
$\|\mathbf{u}\|~= (\mathbf {u},\mathbf {u})^{1/2}$.
  Пространство  Соболева порядка $s\geq 0$ обозначается через
$\mathbf{H}^{s}(G)$, $\|\mathbf {f}\|_s$ -норма его элемента
$\mathbf {f}$;   замыкание в $\mathbf{H}^{s}(G)$ пространства
$\mathcal{C}^{\infty}_0(G)$ обозначается через
$\mathbf{H}^{s}_0(G)$. Пространство Соболева отрицательного порядка
$\mathbf{H}^{-s}(G)$ двойственно к $\mathbf{H}^{s}_0(G)$
 (см. пространство $W_p^{(l)}(\Omega)$ при $p=2$ в $\S 3$ гл. 4 \cite{sob} ,
    $H^k(Q)$ в $\S 4$ гл. 3 \cite{mi}, а также   гл.1 в \cite{mag}).

  В области $G$  с гладкой границей  $\Gamma$ в каждой точке $y\in\Gamma$
  определена нормаль   $\mathbf {n}(y)$  к $\Gamma$.\,
Вектор-функция     $\mathbf {u}$ из $\mathbf{H}^{s+1}(G)$ имеет след
$ \gamma(\mathbf {n}\cdot\mathbf {u})$ на $\Gamma$ ее нормальной
компоненты, который принадлежит пространству Соболева-Слободетцкого
$\mathbf{H}^{s+1/2}(G)$, $|\gamma(\mathbf {n}\cdot\mathbf
{u})|_{s+1/2}$- его норма.

Пусть функция $h\in {H}^{1}(G)$, а $\mathbf{u}=\nabla h$ - ее
градиент. Через ${\mathcal{{A}}}(G)$ обозначим подпространство
$\{\nabla h, h\in H^1(G)\}$ в $\mathbf{L}_{2}(G)$, а через
${\mathcal{{B}}}(G)$ -  его ортогональное дополнение. Так что
\begin{equation}\label{wor  1}\mathbf{L}_{2}(G)=
{\mathcal{{A}}}(G)\oplus{\mathcal{{B}}}(G).
 \footnote {{\small Это разложение я взял из статьи Z.Yoshida и Y.Giga
  \cite{giyo}. Они называют его разложением Вейля \, (H.Weyl\,
  \cite{hw}),
  а пространство ${\mathcal{{B}}}(G)$ обозначают как
  ${L}_{\sigma}^2(G)$. В обобщенном смысле оно формулируется так:
  \[{L}_{\sigma}^2(G)=\{u\in {L}_{}^2(G),
   div\,u=0,\gamma( {n}\cdot{u})=0 \}.\]}}\end{equation} %\newpage
{\small Отметим, что в  разложении Вейля \cite{hw}, Теорема II, роль
$\mathcal{A}$ играет пространство $\mathfrak{G}$-замыкание в норме
${L}_2$ градиентов функций $\psi\in
\Gamma=\mathcal{C}_0^1(G)$.\newline С.Л.Соболев \cite{so54} приводит
разложение $\mathbf{L}_{2}(\Omega)$ в двух случаях, когда $\Omega$
совпадает со всем пространством и случай ограниченной области
$\Omega$,\,гомео-морфной шару. Множества
 вектор-функций вида $\{\nabla\,\varphi\}$ при
 $\varphi\in \mathcal{C}_0^{\infty}(\Omega)$ и
 $\varphi\in \mathcal{C}^{\infty}(\Omega)$, называемых потенциальными,
 обозначаются как $\widetilde{G}_0$ и  $\widetilde{G}_1$, а
 множества вектор-функций вида $\{\mathrm{rot}\,\psi\}$ при
 $\psi\in \mathcal{C}_0^{\infty}(\Omega)$ и
 $\psi\in \mathcal{C}^{\infty}(\Omega)$, называемых соленоидальными,
 обозначаются как $\widetilde{J}_0$ и  $\widetilde{J}_1$, соответственно. %\newline
 ${G}_0$,   ${G}_1$ и ${J}_0$,   ${J}_1$ - их замыкания в норме
   $H={L}_2(\Omega)$. Автор доказывает, что первом случае
   $H=G\oplus\,J$, где $G={G}_0={G}_1$, $J=J_0=J_1$, а во втором:
  $H=G_0\oplus\,I\oplus\,J_0$, где $I=G_1\cap\,J_1$.

    О.А.Ладыженская в $\S 2$ гл.1 книги \cite{lad}
приводит разложение пространства $\mathbf{L}_{2}(\Omega)$ на два
ортогональных подпространства $\mathbf{G}(\Omega)$ и
$\mathbf{J}^\circ(\Omega)$, где $\mathbf{J}^\circ(\Omega)$ есть
замыкание в норме $\mathbf{L}_{2}(\Omega)$ множества $\mathbf{\dot
J}(\Omega)$ бесконечно дифференцируемых финитных в $\Omega$
соленоидаьных векторов, а $\mathbf{G}(\Omega)$ его ортогональное
дополнение в $\mathbf{L}_{2}(\Omega)$. Она пишет: "имеются разные
способы определения этих пространств (см. прежде всего работу Вейля
\cite{hw}, а также \cite{kkr} и \cite{bs})". \,\,Э.Быховский и
Н.Смирнов \cite{bs} отмечают: "К.Фридрихс \cite{fri} доказал близкие
результаты для дифференциальных форм на римановых многообразиях".}

Мы будем придерживаться разложения \eqref{wor  1}.

 Если граница
области $G$ имеет положительный род $\rho$, то $ \mathcal{B}$
содержит в себе конечномерное подпространство
\begin{equation}\label{bh  1} \mathcal{B}_H=\{\mathbf{u}\in\mathbf{L}_{2}(G):
\mathrm{rot}\,\mathbf{u}=0,\,\,\mathrm{div}\,\mathbf{u}=0,
\,\,\gamma(\mathbf{n}\cdot \mathbf{u})=0 \}. \end{equation} Его
размерность
 равна $\rho$ \cite{boso}, а базисные функции
 $h_j\in \mathcal{C}^\infty(G)$. Гладкость обобщенных рещений системы \eqref{bh  1}
 доказали С.Соболев \cite{so37} и Г.Вейль \cite{hw}.   Ортогональное дополнение $\mathcal{B}_H$ в
$\mathcal{B}(G)$, следуя А.Фурсикову \cite{afu}, обозначим через
$\mathbf {V} ^{0} (G)$. Значит,
\begin{equation}\label{bhr  1}{\mathcal{B}(G)=
\mathbf {V} ^{0} (G)\oplus\mathcal{B}_H} (G). \footnote {В
\cite{giyo}:\,
${L}_{\sigma}^2(G)={L}_{\Sigma}^2(G)\oplus{L}_{H}^2(G)$. Символ $L$
перегружен. Мы изменили авторские обозначения пространств
${L}_{\Sigma}^2(G)$
 и ${L}_{H}^2(G)$  на
   $\mathbf {V}^{0} (G)$ и $\mathcal{B}_{H} (G)$. }\end{equation}

   Индексом $\gamma$  будем снабжать пространства вектор-функций
        $\mathbf {u}$, нормальные компоненты которых % из имеет
     имеют на $\Gamma$ нулевой    след
     $ \gamma(\mathbf {n}\cdot\mathbf {u})$:
     \begin{equation}\label{ag  1}\mathcal{A}_{\gamma} =
     \{\nabla h, h\in H^1(G),  \gamma(\mathbf {n}\cdot {\nabla\,h})=0
     \},\end{equation}
     \[ \newline\mathbf{H}^{s+1}_\gamma(G)=\{\mathbf {u}\in \mathbf{H}^{s+1}(G):
     \gamma(\mathbf {n}\cdot\mathbf {u})=0\}.\]
\subsection {Свойства операторов ротор и  градиент дивергенции}
Операторы градиент, ротор и   дивергенция определяются в трехмерном
векторном анализе \cite{zo}. Им соответствует оператор $d$ внешнего
диф-ференцирования на формах $\omega^k$ степени $k=0,1$ и 2.
Соотношения $dd\omega^k=0$ при $k=0,1$ имеют вид
$\mathrm{rot}\,\nabla h=0$ и
$\mathrm{div}\,\mathrm{rot}\,\mathbf{u}=0$.

Формулы $\mathbf{u}\cdot\nabla
h+h\mathrm{div}\mathbf{u}=\mathrm{div}(h \,\mathbf{u}) $,\,
$\mathbf{u}\cdot\mathrm{rot}\,\mathbf{v}-
\mathrm{rot}\,\mathbf{u}\cdot\mathbf{v}=\mathrm{div}[\mathbf{v},\mathbf{u}]$,
где $[\mathbf{v},\mathbf{u}]$ - векторное произведение, и
интегрирование по области $G$ используются при определении
операторов $\mathrm{div}\,\mathbf{u}$  и $\mathbf{rot}\mathbf{u}$ в
$\mathbf{L}_{2}(G)$.

 Оператор Лапласа
выражается через $\mathrm{rot}\,\mathrm{rot}$ и
$\nabla\,\mathrm{div}$:
\begin{equation}\label{dd  1}-\Delta \mathbf{v} =\mathrm{rot}\, \mathrm{rot}\,\mathbf{v}-
  \nabla \mathrm{div} \mathbf{v}.\end{equation}
  Оператор Лапласа эллиптичен, а операторы $\mathrm{rot}$ и
  $\nabla\,\mathrm{div}$ не являются таковыми \cite{so71}. Они
вырождены, причем  $\mathrm rot\,\mathbf {u}=0$ при $\mathbf{u}\in~
\mathcal{A}(G)$, $\nabla\mathrm div\,\mathbf {v}=0$ при
$\mathbf{v}\in \mathcal{B}(G)$ в смысле теории распределений.

Операторы $\mathrm{rot} +\lambda\, I$ на $\mathcal{A}(G)$ и
$\nabla\mathrm{div} +\lambda\, I$ на $\mathcal{B}(G)$ сводятся к
  умножению на $\lambda\neq 0$, а на ортогональных подпространствах
   при обращении  требуют краевых условий; например,
    $ \gamma(\mathbf {n}\cdot\mathbf {u})=0$.
 Они принадлежат классу Б.Вайберга и В.Грушина \cite{vagr}
операторов, приводимых  к эллиптическим, так как их расширения
являются эллиптическими переопределенными операторами. Краевые
задачи с условием
 $ \gamma(\mathbf {n}\cdot\mathbf {u})=g$ являются эллиптическими по
 Солонникову  \cite{so71}.% (см. пп.2.1 и 4.1).
  В пространства $\mathbf{V}^{0}(G)$   и
${\mathcal{{A}}_{\gamma}}(G)$ операторы ротор и градиент дивергенции
допускают самосопряженные расширения.

А именно, оператор \quad $S:\mathbf{V}^{0}(G)\longrightarrow
\mathbf{V}^{0}(G)$ с областью определения \quad
$\mathbf{W}^{1}=\{\mathbf {u}\in\mathbf {V}^0(G):\mathrm{rot}\mathbf
{u}\in\mathbf {V}^0(G)\}$,  совпадающий с $\mathrm{rot}\,\mathbf{u}$
при  $\mathbf {u}\in\mathbf{W}^{1}$, %$\subset\mathbf{H}^{1}(G)$,
 является самосопряженным и имеет вполне непрерывный обратный  $S^{-1}$ из
$\mathbf{V}^{0}$ в $\mathbf{W}^{1}$ \cite{giyo}.

Соответственно, оператор
 $\mathcal{N}_d:\mathcal{A}_{\gamma}(G)
 \longrightarrow \mathcal{A}_{\gamma}(G)$ , совпадающий с
 $\nabla \mathrm{div}\mathbf {u}$
 на  $\mathcal{A}^{2}=\{\mathbf {u}\in\mathcal{A} _ {\gamma}(G):
 \nabla \mathrm{div}\mathbf {u}\in\mathcal{A} _ {\gamma}(G)\}$,  самосопряжен и
 его обратный оператор
 $\mathcal{N}^{-1}_d:\mathcal{A} _ {\gamma}\rightarrow \mathcal{A}^{2}$
  вполне непрерывен ( п.4.3). %теорема 5

Следовательно, каждый из этих операторов имеет полную систему
собственных функций, отвечающих ненулевым собственным значениям:
\[\mathrm{curl}\mathbf{u}_{j}^{\pm}={\pm}\lambda_j\mathbf{u}_{j}^{\pm},
\quad {\lambda}_j\in {\Lambda}\subset\mathbb{R},\quad \nabla
\mathrm{div}\mathbf{q}_{j}=\mu_j \mathbf{q}_{J},\quad \mu_j\in M
\subset \mathbb{R},\]
\begin{equation}\label{gd  1}\mathbf{a}(x)=\sum_{\mu_j\in M}(\mathbf{a},\mathbf{q}_{j})\mathbf{q}_{j}
(x),\quad \text{если}\quad \mathbf{a}(x)\in{\mathcal {{A}} _
{\gamma}} (G),\quad \|\mathbf{q}_{j}\|=1,\end{equation}
\[\mathbf{b}(x)=\sum_{\lambda_j\in
\Lambda}[(\mathbf{b},\mathbf{u}_{j}^{+})
\mathbf{u}_{j}^{+}(x)+(\mathbf{b},\mathbf{u}_{j}^{-})
\mathbf{u}_{j}^{-}(x)], \quad
%\text{при}\,\,\,
\mathbf{b}(x)\in\mathbf {V} ^ {0} (G),\quad
\|\mathbf{u}_{j}^{\pm}\|=1.\]

 В  шаре $B$ радиуса $R$ собственные функции
 $\mathbf{u}^{\pm}_{\kappa}$ ротора, отвечающие
  ненулевым собственным значениям
  $\pm\lambda_{\kappa}=\pm\rho_{n,m}/R$ и собственные функции
  $\mathbf{q}_{\kappa}$
градиента дивергенции с собственными значениями $\nu_{\kappa}^2$,
$\nu_{\kappa}=\alpha_{n,m}/R,$ выражаются  явными формулами (см.
  пп. 3.2 и 5.3) и %\cite{saUMJ13}
    \[\mathbf{rot}\,\mathbf{u}^{\pm}_{\kappa}=\pm\lambda_{\kappa}\,
\mathbf{u}^{\pm}_{\kappa}, \quad
\gamma\mathbf{n}\cdot\mathbf{u}^{\pm}_{\kappa}=0;\quad
\mathbf{rot}\,\mathbf{q}_{\kappa}=0
 \quad \kappa=(n,m,k),\,\]
\[ \nabla\,div\,\mathbf{u}^{\pm}_{\kappa}=0, \quad
\nabla\,div\,\mathbf{q}_{\kappa}=\nu_{\kappa}^2\mathbf{q}_{\kappa},
\quad \gamma\mathbf{n}\cdot\mathbf{q}_{\kappa}=0, \quad
 \quad\,\, |k|\leq n ,\] где
числа $\pm\rho_{n,m}$ и $\alpha_{n,m}$ - нули функций $\psi_n$ и их
производных $\psi_n'$, а
 \begin{equation}
\label{psi__1_}\psi_n(z)=(-z)^n\left(\frac{d}{zdz}\right)^n\frac{\sin
z}z, \quad n\geq 0,\,\,m\in \mathbb{N}.\end{equation}

Cобственные функции каждого из операторов взаимно ортогональны и их
совокупная система полна в ${\mathbf{{L}}_{2}}(B)$ \cite{saUMJ13}.

Найдены необходимые и достаточные условия на вектор-функции
$\mathbf{u}$ из $\mathbf{V}^{0}(B)$ и $\mathbf{v}$ из
$\mathcal{A}_{\gamma}(B)$, при которых их ряды Фурье сходятся в
норме пространства Соболева $\mathbf{H}^{s}(B)$ порядка $s>0$.
  Они состоят в  принадлежности $\mathbf{u}$ и $\mathbf{v}$
пространствам $\mathbf{V}^{s}_{\mathcal{R}}(B)$ п. 3.3 и
$\mathcal{A}_{\mathcal{K}}^{s}(B)$ п. 5.7.

  Предлагается единый подход к изучению краевых задач
  \eqref{kron__1_},\eqref{ndi__1_}
в пространствах $\mathbf{H}^{s}(G)$, $s\geq 1$, при $\lambda\neq 0$.
 Их разрешимость  зависит от пространств, к
которым принадлежат $\mathbf{f}$ и $g$ ( пп.2.3 и 4.3 ).

 \subsection{Спектральная задача} Пусть  $G$ - ограниченная область
  в ${{R}^{3}}$ с гладкой
границей $\Gamma $, $\mathbf{n}$- внешняя нормаль к $\Gamma $. В
частности, $G$ может быть шаром $B$,  $|x|<R$,
 с границей $S$.
\begin{problem}
Найти все
собственные значения $\lambda $ и собственные вектор-функции
$\mathbf{u}(\mathbf{x})$ в ${{\mathbf{L}}_{2}}(G)$ оператора ротор
такие, что
 \begin{equation}
\label{ro__1_} \mathrm{rot}\,\mathbf{u}=\lambda \mathbf{u}\quad
\text{в}\quad G, \quad \mathbf{n}\cdot \mathbf{u}|_{\Gamma }=0,
\end{equation}
 где $\mathbf{n}\cdot \mathbf{u}$ - скалярное произведение векторов
$\mathbf{u}$ и $\mathbf{n}$.\end{problem}

 К области определения $\mathcal{M}_{\mathcal{R}}$
 оператора $\mathcal{R} $ задачи 1
    отнесем все вектор-функции
 $\mathbf{v}(\mathbf{x})$ класса $\mathcal{C}^2(G)\cap
 \mathcal{C}(\overline{G})$, удовлетворяющие граничному условию
  и такие, что $\text{rot}\,\mathbf{v}\in
 {\mathbf{L}}_{2}(G)$.
  Пространство основных вектор-функций ${\mathcal{D}}(G)$
 содержится в $\mathcal{M}_{\mathcal{R}}$ и
 плотно в ${{\mathbf{L}}_{2}}(G)$ \cite{vla}. % ${\mathbf{D}}(G)$
     \subsection{О приложениях}  Собственные функции задачи 1
     имеют приложения в
       гидродинамике,
        где они называются полями Бельтрами; % \,\cite{lad},\cite{koz};
         в астрофизике         и в физике плазмы
         они называются бессиловыми полями (force-free magnetic
         fields -
           С. Чандрасекхар и П.Кендал  \cite{chake}, free-decay fields -
                            Д. Тэйлор \cite{tay}).
 В теоремах В.И.Арнольда \cite{ar}\, 1965, и В.В.Козлова, 1983,
 (см. \cite{koz}) изучавших
  топологию линий тока течений
 идеальной жидкости, имеется условие
  $[\text{rot }\mathbf{v}, \mathbf{v}]\neq 0$.
   Стационарные течения вязкой несжимаемой  жидкости со    скоростью
$\mathbf{v}(\mathbf{x})$, удовлетворяющей уравнению
 $\text{rot }\mathbf{v}= \lambda\,\mathbf{v}$
  (в классе периодических функций)
 изучал М. Энон \cite{he}.
  Ссылаясь на его вычисления,
   В.~ Арнольд пишет, что такие течения
   "могут иметь линии тока с весьма сложной топологией,
   характерной для задач небесной механики".

Самосопряженные расширения оператора ротор и его свойства %задачи 1
изучали П.Е.Берхин \cite{berh}\, 1975, \,  З.~Иошида и И.~Гига
\cite{giyo}\, 1990
 а также Р.~Пикар \cite{pi}\, 1996,  и  Н.Филонов \cite{phil}.

 Д. Кантарелла
и Де Турк, Г.Глюк и М.Тэйтель 2000 \cite{cdtgt} исследовали
топологию линий тока собственных функций ротора с минимальным
собственным значением в шаре и в шаровом слое.

 С.Чандрасекхара и П.Кендала \cite{chake}\, 1957, заметили, что
 собственные  функции ротора  можно выразить через решения уравнения
 Гельмгольца. В цилиндре ( с условием периодичности
вдоль оси), эта идея была реализована в работе Д. Монтгомери,
Л.Тернера и Г.Вахалы о магнито-гидродинамической турбулентности
\cite{motuva} \,1978. Они отмечают: " три интегральных инварианта (
полная энергия, магнитная   и крос спиральность) имеют простые
квадратичные выражения в терминах коэффициентов разложения в ряды
Фурье".

 Другие приложения собственных функций  и рядов Фурье оператора ротор
 имеются в работах автора.\,{\small В частности, найдена  связь
между собственными функциями операторов ротора и Стокса, построены
явные решения нелинейных уравнений Навье-Стокса, разработан метод
численного решения задачи Коши для уравнений Навье-Стокса
\cite{sa04}--\cite{saUMJ13}.}{\footnote {\small В 2003 году О.А.
Ладыженская решала задачу "О построении базисов в пространствах
соленоидальных векторных полей"\,\cite{lab}  и искала способы
вычисления собственных функций оператора Стокса в областях
простейших форм  в явном виде. Автору удалось найти их в случае
периодических граничных условий и в шаре \cite{sa04},\cite{sa07}. }}
 \subsection{ Структура работы и основные результаты}
В $\S\, 2$ в ограниченной области $G$ с гладкой границей $\Gamma$
исследована краевая задача
\begin{equation}
\label{kron__1_}\text{rot}\mathbf{u}+\lambda\,
\mathbf{u}=\mathbf{f}(\mathbf{x}), \quad
 \mathbf{x}\in G, \quad \mathbf{n}\cdot \mathbf{u}{{|}_{\Gamma }}=g
 \end{equation}
 для
 оператора ротор в пространствах $\mathbf{H}^{s+1}(G)$, число $s+1\geq 1$ целое.
Определяется ее оператор $\mathbb{A}$\,(см.\eqref{op  1}).

 Доказано, что при $\lambda\neq
0$ эта задача  является обобщенно эллиптической: она приводится к
 переопределенной эллиптической задаче по определению Солонникова
\cite{so71}. Из его Теоремы 1.1 вытекает Теорема 1 и, в частности,
конечномерность ядра оператора $\mathbb{A}$ задачи в пространстве
Соболева $\mathbf{H}^{s+1}(G)$ и априорная оценка:
\begin{equation}
\label{arot__1_} C_s\|\mathbf{u}\|_{s+1}
\leq\|\mathrm{rot}\,\mathbf{u}\|_{s}+
|\lambda|\|\mathrm{div}\,\mathbf{u}\|_{s}+
|\gamma({\mathbf{n}}\cdot\mathbf{u})|_{s+1/2}+
\|\mathbf{u}\|_{s}.\end{equation}
 В п.2.3 мы изучаем оператор
$\mathrm{rot}~+\lambda \mathbf{I} $  в ортогональных
подпространст-вах в $\mathbf{L}_2(G)$. На %подпространствах
$\mathcal{A}$ и $B_H$ он    сводится к
   $\lambda\,\mathbf{u} $.
   На $\mathbf{V}^{0}(G)$ он
продолжается как самосопряженный оператор
$S+\lambda\,I$
%$:\mathbf{V}^{0}(G)\rightarrow\mathbf{V}^{0}(G)$
%\cite{giyo}.
 Выписаны необходимые и достаточные условия его
обратимости \, (Теорема 2).

В $\S\, 3$ указывается способ решения спектральной задачи 1  в шаре.
При $\lambda\neq 0$ задача сводится к спектральной задаче Дирихле
для скалярного оператора Лапласа с условием $v(0)=0$ в  центре шара.
Она решается явно \cite{vla}, это позволяет определить радиальные
компоненты собственных вектор-функций и числа $\lambda^2_\kappa>0$.
Две другие компоненты определяются из уравнений
$\mathrm{rot}\mathbf{u}_\kappa=\pm\lambda_\kappa\mathbf{u}_\kappa$,
$\mathrm{div}\mathbf{u}_\kappa=0$\, .
 Ее решение опубликовано\,\cite{saUMJ13}, в
  п.3.2 мы приводим уравнения для ненулевых собственных значений
   и  формулы собственных функций.

  В.П.Михайлов \cite{mi} выделил подпространства
   $H^s_\mathcal{D}(G)$ и $ H^s_\mathcal{N}(G)$ в пространстве Соболева
$ H^s(G)$ и доказал, что условие принадлежно-сти $f$ к $
H^s_\mathcal{D}(G)$
 (соотв., к $ H^s_\mathcal{N}(G)$)   необходимо и достаточно для
 сходимости ее ряда Фурье    по системе собственных функций
    оператора    Лапласа c условием Дирихле ( Неймана)
     в норме $ H^s(G)$.

 Мы приводим в п.3.3
 аналогичный результат для оператора ротор (Теорема 3 ), а  п.5.6
 -- для  градиента дивергенции.

В $\S\, 4$ исследована краевая задача
\begin{equation}
\label{ndi__1_}\nabla\text{div}\mathbf{u}+\lambda\,
\mathbf{u}=\mathbf{f}(\mathbf{x}), \quad
 \mathbf{x}\in G, \quad \mathbf{n}\cdot \mathbf{u}{{|}_{\Gamma }}=g
 \end{equation}
    в
ограниченной области $G$ с гладкой границей $\Gamma$. Доказано,что
эта задача обобщенно эллиптична при $\lambda\neq 0$: она приводится
к эллиптической задаче \cite{so71}. Откуда вытекает Теорема 4,
 конечно-мерность ядра оператора $\mathbb{B}$ задачи в
пространстве Соболева $\mathbf{H}^{s+2}(G)$ и априорная оценка:
%\eqref{ond__2_}.
\begin{equation}
\label{ond__1_} C_s\|\mathbf{u}\|_{s+2}
\leq|\lambda|\|\mathrm{rot}^2\,\mathbf{u}\|_{s}+
\|\nabla\mathrm{div}\,\mathbf{u}\|_{s}+
|\gamma({\mathbf{n}}\cdot\mathbf{u})|_{s+3/2}+ \|\mathbf{u}\|_{s}.
\end{equation}
В п.4.3 мы изучаем оператор $\nabla\mathrm{div} +\lambda \mathbf{I}
$  в ортогональных подпространствах $\mathcal{A}$ и $\mathcal{B}$ в
$\mathbf{L}_2(G)$. На  $\mathcal{B}$ %подпространстве
  оператор $\nabla\mathrm{div} \mathbf{u}
+\lambda \mathbf{u} $ сводится к  %оператору
 $\lambda\mathbf{u}$. На
 $\mathcal{A}_{\gamma}$  он продолжается
 как самосопряженный оператор $\mathcal{N}_d + \lambda\,I$
 %$: \mathcal{A}_{\gamma}\rightarrow\mathcal{A}_{\gamma}$.
 Найдены необходимые и достаточные условия его
обратимости (Теорема 5b).

В  $\S\, 5$ спектральная задача для оператора градиент дивергенции в
  области с гладкой границей сводится к решению спектральной
задачи Неймана для скалярного оператора Лапласа.

 В шаре ее решения
вычислены явно\,\cite{vla}. В результате мы получаем формулы
\eqref{qom 1} собственных функций $\mathbf{q}_{\mu}(\mathbf{x})$
 градиента дивергенции.

 В   $\S 6$ мы рассматриваем совокупную систему
  собственных функций  ротора и градиента дивергенции:\quad
$\{\mathbf{q}_{i}(\mathbf{x}), \,\mathbf{h}_{j}(\mathbf{x}),\,
 \mathbf{q}_{k}^{+}(\mathbf{x}),\, \mathbf{q}_{k}^{-}
 (\mathbf{x})\}$ \footnote {{\small Вектор-функции
$\mathbf{q}_{i}(\mathbf{x})$ удовлетворяют также уравнениям
$\bf{rot}\mathbf{q}_{i}=0$, а $\mathbf{q}_{k}^{\pm}(\mathbf{x})$ -
уравнениям $\nabla div\mathbf{q}_{k}^{\pm}(\mathbf{x})=0$. }}
\newline $\mu_i\in M,\,\, j\in [1,\rho],\, \,\lambda_k\in \Lambda,$
они взаимно ортогональны и образует в  ${\mathbf{{L}}_{2}}(G)$
ортонормированный базис.

 В шаре $B$ векторное поле $\mathbf{f}\in \mathbf{L}_{2}(B)$ разлагается на
 потенциальное и соленоидальное поле $\mathbf{a}_f$ и $\mathbf{b}_f$:
 $\mathbf{f}(\mathbf{x})=\mathbf{a}_f(\mathbf{x})+
 \mathbf{b}_f(\mathbf{x})$.% \cite{lad, bs}.

В качестве примера в $\S\, 7$   методом Фурье решена краевая задача
\eqref{grd__1_}  в шаре при любых $\lambda$ и $g=0$   (Теорема 8).

\section{ Оператор ротор в ограниченной области}

\subsection{Краевая задача:} в ограниченной области $G$ с гладкой границей
$\Gamma$ заданы векторная и скалярная функции $\mathbf{f}$ и ${g}$,
найти вектор-функцию $\mathbf{u}$, такую что
\begin{equation}
\label{grd__1_}
  \text{rot}\mathbf{u}+\lambda\, \mathbf{u}=\mathbf{f}(\mathbf{x}),   \quad
 \mathbf{x}\in G, \quad \mathbf{n}\cdot \mathbf{u}{{|}_{\Gamma }}=g.
\end{equation} Эта задача не эллиптична.
Оператор $\text{rot }+\lambda\text{I}$ первого порядка не является
эллиптическим, так как ранг его символической матрицы
$\mathrm{rot}(i\xi)$,  равный двум при всех $\xi\in
\mathcal{R}^3\backslash 0$, меньше трех \cite{sa75}.

 Б.Вайнберг и В.Грушин   \cite{vagr} 1967
определили   на гладком многообразии $X$ без края  класс {\it
равномерно неэллиптических систем}  (РНС) сингуляр-ных
интегро-дифференциальных уравнений и класс матричных с.и.д.
операторов, {\it глобально приводимых к эллиптическим матрицам}, и
доказали их эквивалентность. Эти определения требуют введения
дополнительных понятий.

 Мы приведем их  для
систем дифференциальных уравнений с постоянными коэффициентами,
который обозначим как (РНСp)

 Система дифференциальных уравнений, $S(D)u=f$ порядка $m$, из этого
класса обладает свойствами: \newline
 а)\,ее символическая матрица
$S_0(i\xi)$ имеет постоянный ранг при всех
$\xi\in\mathcal{R}^3\backslash 0$.

 Это позволяет построить
аннулятор $C(D)$ оператора $S_0(D)$ такой, что $(CS_0)(D)\equiv 0$
на $X$ и определить
\newline б)\, расширенную систему \quad $Su=f,// CSu=Cf$ порядка~$m$.
\newline
Ее символическая матрица $S_0(i\xi),//(CS)_0(i\xi)$ определяется
младшей частью оператора $S(D)$ и дополняет матрицу $S_0(i\xi)$.
\newline в)\,Если ранг  расширенной матрицы максимален, то
исходная система $Su=f$ принадлежит классу (РНС1) и  степень ее
приводимости равна единице.

г)\, Если система $Su=f$ такова, что ранг  расширенной матрицы не
максимален, но постоянный, то процесс повторяется и при определенных
условиях система принадлежит классу (РНС2). И  так далее.

Авторы \cite{vagr} доказали, что система $Su=f$  класса (РНСp)
являются разрешимой по Фредгольму или Нетеру в  пространствах
Соболева $\mathbf{H}^s(X)$, если $f\in \mathbf{H}^{s-m+p}(X)$, где
$s\geq m$ целое. В качестве примера оператора из класса (РНС1) они
приводят оператор $d+\ast$ на дифферен-циальных формах степени $k$ в
$2k+1$-мерном многообразии $X$.\footnote {Другие классы обобщенно
эллиптических операторов см. в работе \cite{sak97}.}

 Покажем, что дифференциальный оператор $(\text{rot}+\lambda\,I)$
   при $\lambda\neq 0$ принадлежит классу (РНС1) в любой
области $X\subset \mathbf{E}^3$. Действительно,

а) его символическая матрица $\mathrm{rot}(i\xi)$ не зависит от $x$
и ее ранг равен двум при всех $\xi\in \mathcal{R}^3\backslash
0$.\newline
 б) оператор $\mathrm{rot} $ имеет левый
аннулятор $\mathrm{div}$ : $\mathrm{div}\,\mathrm{rot}\,
\mathbf{u}=0$ на $X$.
\newline в) ранг символичекой $(4\times 3)$-матрицы
$[\mathrm{rot}(i\xi)// \lambda\,\mathrm{div}(i\xi)]$ равен трем при
всех $\xi\in \mathcal{R}^3\backslash 0$. Следовательно расширенная
система
  \begin{equation} \label{rodi__1_}
\mathrm{rot}\,\mathbf{u}+\lambda \mathbf{u}=\mathbf{f},\quad \lambda
\mathrm{ div}\,\mathbf{ u}=\mathrm{div}\,\mathbf{f},\end{equation}
является эллиптической системой первого порядка, а система
\eqref{grd__1_} принадлежит классу (РНС1).

Далее,  система  \eqref{rodi__1_} с произвольной функцией $f_4$
вместо $\mathrm{div}\,\mathbf{f}$ и с краевым условием
$\gamma\,\mathbf{n}\cdot \mathbf{u}=g$
  составляют   переопределенную эллиптическую краевую задачу
 по Солонникову \cite{so71}. А именно,

1) система  \eqref{rodi__1_} эллиптична,

2) оператор краевого условия $\gamma\mathbf{n}\cdot \mathbf{u}$\,
"накрывает" \,\,оператор системы. % \eqref{rodi__1_}.

Первое условие сводится к тому, что  однородная система линейных
алгебраических уравнений:
\begin{equation}
\label{cdx  1}\mathrm{rot}(i\xi )\mathbf{w}=0, \quad
\lambda\,\mathrm{div}(i\xi )\mathbf{w}=0, \quad \forall \xi\neq 0
\end{equation}
 c параметром $\xi \in \mathcal{R}^3$ имеет только тривиальное
 решение $ \mathbf{w}=0$.

  Второе условие означает, что  однородная система линейных
диф-ференциальных уравнений:
\begin{equation}\label{cdz  2}
\mathrm{rot}(i\tau+\mathbf{n} d/dz ) \mathbf{v}=0,\quad
\mathrm{div}(i\tau+\mathbf{n} d/dz )\mathbf{v}=0, \quad \forall \tau
\neq 0,\end{equation} на полуоси $z\geq 0$ с краевым условием:
$\mathbf{n}\cdot \mathbf{v}|_{ z=0}=0$ и убыванием,
$\mathbf{v}(y,\tau; z )\rightarrow 0$ при $z\rightarrow + \infty$,
имеет только тривиальное решение.

Здесь  $\tau$ и $\mathbf{n}$ -- касательный и нормальный векторы к
 $\Gamma$ в точке $y\in \Gamma$ и $|\mathbf{n}|=1$.
 \footnote {  Главные части  системы в  \cite{so71}
 определяются с помощью весов $s_k$ и $t_j$ таких, что
 $ord\, L_{k,j}\leq s_k+t_j $.
 %матрицы $ L $ не превосходит   $$. \quad
    Положив    $s_k=0$\, при $k=1-4$
    и $t_j~=1$ при $j=1-3$ мы получим операторы
    системы  \eqref{cdx  1}, а в краевом операторе -
    %нулевого    порядка полагаем
     $\sigma_1=-1$.}

При доказательстве утверждений  1),2) воспользуемся соотношением
\begin{equation}\label{cc  1} \mathrm{rot}\, \mathrm{rot}\,\mathbf{v}= -\Delta \mathbf{v}
 + \nabla \mathrm{div} \mathbf{v}.\end{equation}
Тогда $1^0)$. Из уравнений  (\ref{cdx  1}) вытекает уравнение
$-\Delta(i\xi)\mathbf{w}=0$. Оно распадается на три скалярных
уравнения $|\xi | ^2 w_j =0$.% \newline
Значит, $\mathbf{w}=0$ при $|\xi | \neq 0$. Эллиптичность системы
(\ref{rodi__1_}) доказана.
\newline
$2^0)$. Из уравнений (\ref{cdz  2}) получаем уравнение $ (-|\tau| ^2
+ (d/dz)^2)\mathbf{v} = 0$ с параметром $|\tau |> 0$. Его убывающее
при $z\rightarrow + \infty$ решение имеет вид:
$\mathbf{v}=\mathbf{w} e^{-|\tau|z}$ . Оно удовлетворяет уравнениям
(\ref{cdz 2}), если вектор-функция $\mathbf{w}$ такова, что
  $ \omega\times \mathbf{w}=0,\quad {\omega}' \cdot
\mathbf{w}=0$,   где $ \omega\equiv i\tau -| \tau| \mathbf{n}$
--вектор-столбец,  $ \omega'$ -- вектор-строка, а ${\omega}' \cdot
{\omega}$ -- их произведение.

Легко убедиться, что векторное и  скалярное произведения $\omega$
 на $\omega$ равны нулю: $\omega \times \omega=0,$ \, $\omega'\cdot\omega=0$.
Ранг матрицы $\mathrm{rot} (i\xi )$ равен двум при  $\xi \neq 0$,
поэтому  $\mathbf{w}=c\,\omega$, где $c$ - постоянная,  и других
решений нет.
 Граничное условие приводит нас к уравнению:
  $|\tau |c=0$.   Значит  $c=0$ при $|\tau| > 0$ и, следовательно, $\mathbf{v}=0$.

Итак, система  \eqref{rodi__1_} с краевым условием $\mathbf{n}\cdot
\mathbf{u}{{|}_{\Gamma }}=g$  при $\lambda\neq 0$ является
эллиптической задачей.

Мы скажем в этом случае, что задача \eqref{grd__1_}  при
$\lambda\neq 0$ является обобщенно эллиптической.
%%%%2.2
\subsection{Оператор задачи в пространствах Соболева}
Пусть вектор-функция $\mathbf{u}$ принадлежит пространству Соболева
${\bf{H}^{s+1}}(G)$, где $s\geq 0$--целое. Тогда компоненты
$\mathrm{rot} \mathbf{u}$ и $\mathrm{div} \mathbf{u}$ принадлежат
${H}^{s}(G)$, а вектор-функция $\mathbf{f}:=\mathrm{rot} \mathbf{u}+
\lambda \mathbf{u}$ принадлежит пространству
\begin{equation}\label{pr  2} {\bf{E}^{s}}(G)=\{\mathbf{f}\in
{\mathbf{H}^{s}(G)}: \mathrm{div}\,\mathbf{f}\in
 {H}^{s}(G)\},\end{equation} которое снабжается нормой
 $\|\mathbf{v}\|_{\mathbf{E}^{s}}=(
 \|\mathbf{v}\|^2_{s }+
 \|\mathrm{div}\mathbf{v}\|^2_{s})^{1/2}$.

Далее $g:=\gamma({\mathbf{n}}\cdot\mathbf{u})\equiv\mathbf{n}
\cdot\mathbf{u}|_{\Gamma}$ принадлежит пространству
Соболева-Слободетского $H^{s+1/2}(\Gamma)$. \newline Следовательно,
при $\lambda\neq 0$ задаче соответствует ограниченный оператор
\begin{equation}\label{op  1} \mathbb{A}\mathbf{u}\equiv \begin{matrix}
\mathrm{rot}\,\mathbf{u}+\lambda\,\mathbf{u} \\\gamma
(\mathbf{n}\cdot \mathbf{u})\end{matrix}: \bf{H}^{s+1}(G)\rightarrow
\begin{matrix}\bf{E}^{s}(G)\\ H^{s+1/2}(\Gamma)\end{matrix}.
\end{equation}
Согласно Теореме 1.1 из работы Солонникова \cite{so71} о
переопределенных эллиптических краевых задачах в ограниченной
области $G$ с гладкой границей $\Gamma\in \mathcal{C}^{s+1} $,
обобщенно эллиптический оператор \eqref{op 1} имеет левый
регуляризатор: то-есть
 ограниченный оператор  $\mathbb{A}^L$ такой, что
$\mathbb{A}^L\mathbb{A}=\mathbb{I}+\mathbb{T}$, где  $\mathbb{I}$ -
единичный, а $\mathbb{T}$ - вполне непрерывный операторы, и
существует постоянная $C_s >0$ такая, что выполняется априорная
оценка:
\begin{equation} \label{apro s} C_s\|\mathbf{u}\|_{s+1}
\leq\|\mathrm{rot}\,\mathbf{u}\|_{s}+
|\lambda|\|\mathrm{div}\,\mathbf{u}\|_{s}+
|\gamma({\mathbf{n}}\cdot\mathbf{u})|_{s+1/2}+ \|\mathbf{u}\|_{s}.
 \end{equation}
Оценка \eqref{apro s} известна (см. например \cite{bobr,giyo}). Мы
показали, что она получается из работы В.А.Солонникова \cite{so71}.
{\footnote  {\small Он приводит  оценку в банаховых пространствах
Соболева ${{W}^{l+t_j}_p}(G)$, $l\geq 0, p>1$, которые при $l=s,
p=~2,t_j=~1$ совпадают %гильбертовыми пространствами
${{H}^{s+1}}(G)$, и доказывает, что эта оценка является точной.}}
%\cite{so71}, \cite{rt},\cite{sa75}${\bf{H}^{s+1}}(G)$,

 Линейное пространство  решений  однородной задачи
  обозначим через $\mathcal{N}$.\quad
 Итак, имеет место
\begin{theorem}
Оператор $\mathbb{A}$ в пространствах (\ref{op  1}) имеет левый
регуля-ризатор.
%\newline
 Его ядро $\mathcal{N}$ конечномерно и
выполняется  оценка \eqref{apro s}.
\end{theorem}

Из этой теоремы и  оценки следует, что при $\lambda \ne 0$

 a){\it число линейно независимых решений %спектральной
задачи 1 конечно,}

 b){\it любое (обобщенное) решение задачи бесконечно дифференцируемо вплоть до
границы, если граница области бесконечно дифференцируема.}
\subsection{Оператор $\mathrm{rot}  +\lambda
\mathbf{I} $ в подпространствах $\mathbf{L}_2(G)$}
На подпространстве $\mathcal{A}$ %=\{\mathbf{u}=\nabla\,h:\, h\inH^1(G)\}$
 оператор $\mathrm{rot} \mathbf{u} +\lambda \mathbf{u} $
сводится к $\lambda \mathbf{u} $.

Ортогональное дополнение $\mathcal{B}$ к %подпространству
$\mathcal{A}$ в %пространстве
$\mathbf{L}_{2}(G)$ определяется так
\begin{equation}\label{oB 1}
 \mathcal{B}=\{
\mathbf{u}\in
\mathbf{L}_{2}(G):\int\limits_{G}\mathbf{u}\cdot\nabla\,h\,d\,x=0,\quad
\text{для любой}\quad h\in H^1(G)\}.\end{equation}
  Для функций $\mathbf{u}$ из $\mathbf{H}^1(G)$ получаем:  $\mathrm{div}\,\mathbf{u}=0$ в $G$
 и $\mathbf{n}\cdot
\mathbf{u}|_{\Gamma}=0$.

В пространстве $\mathcal{B}$ выделяется подпространство
\begin{equation}\label{oBH 1}
 \mathcal{B}_H=\{\mathbf{u}\in\mathcal{B}:
\int\limits_{G}\mathbf{u}\cdot\mathrm{rot}\,\mathbf{v}\,
d\,\mathbf{x}=0,\quad \text{для любой}\quad \mathbf{v}\in
\mathcal{D}(G)\}.\end{equation}
 Пространства $\mathcal{B}$ и $\mathcal{B}_H$ в обобщенном смысле
 обозначают  так:
 \[ \mathcal{B}=\{\mathbf{u}\in\mathbf{L}_{2}(G):
\mathrm{div}\,\mathbf{u}=0, \,\,\mathbf{n}\cdot
\mathbf{u}|_{\Gamma}=0\},\]
 \begin{equation}\label{oBH 2}
 \mathcal{B}_H=\{\mathbf{u}\in\mathbf{L}_{2}(G):
\mathrm{rot}\,\mathbf{u}=0,\,\,\mathrm{div}\,\mathbf{u}=0,
\,\,\mathbf{n}\cdot \mathbf{u}|_{\Gamma}=0\}.
 \end{equation}
%Пространство  конечномерно.
Ввиду оценки
 \eqref{apro s}  базис $\mathcal{B}_H$ состоит из
 бесконечно дифференцируемых в $G$ вектор-функций $\{\mathbf{h}_j\}_
 {j\in\overline{1\rho}}$, где $\rho$ есть род границы $G$ \cite{ boso}.

  Ортогональное дополнение к
 $\mathcal{B}_H$ в $\mathcal{B}$ обозначим как $\mathbf{V}^0(G)$,
причем
$\|\mathbf{u}\|_{\mathbf{V}^0}=\|\mathbf{u}\|_{\mathbf{L}_2}$. Так,
что
\begin{equation}\label{cB 0}
\mathbf{L}_{2}(G)=\mathcal{A}\oplus\mathcal{B}, \quad
\mathcal{B}=\mathcal{B}_H\oplus\mathbf{V}^0(G) .\end{equation} В
случае шара % пространство $\mathcal{B}_H$ пусто и
$\mathcal{B}=\mathbf{V}^0(G)$. В полнотории
 $dim\,\mathcal{B}_H =1$.

  Наконец, в $\mathbf{V}^0(G)$ выделяется
подпространство
\begin{equation}\label{pW  1}
 \mathbf{W}^1(G)=\{
\mathbf{u}\in \mathbf{V}^{0}(G):\mathrm{rot}\,\mathbf{u}\in
\mathbf{V}^{0}(G) \}.\end{equation}
 В  силу
 оценки \eqref{apro s} оно содержится $\mathbf{H}^1(G)$  и плотно
 в $\mathbf{V}^{0}(G)$,
так как плотное в нем множество
 $\mathbf{C}_0^{\infty}\cap\mathbf{V}^{0}(G)$
содержится в $\mathbf{W}^1(G)$.

 З.~Иошида и И.~Гига  \cite{giyo} определили в гильбертовом пространстве
$\mathbf{V}^{0}(G)$ оператор $S:\mathbf{V}^{0}(G)\rightarrow
\mathbf{V}^{0}(G)$, который совпадает с $\mathrm{rot}\,\mathbf{u}$
при $\mathbf{u}\in\mathbf{W}^1(G)$, и доказали, что

 {\it  Оператор $S$ является
 самосопряженным и  имеет вполне непрерывный обратный оператор $S^{-1}$ из
  $\mathbf{V}^0(G)$   в $\mathbf{W}^1(G)$. Спектр  $\sigma(S^{-1})$
   точечный и действительный
   и не содержит точек накопления кроме нуля.
 Семейство собственных функций оператора
  $S$ образует
  ортогональный базис в  пространстве
   $\mathbf{V}^0(G)$. }

Собственные функции оператора
  $S$ принадлежат пространствам  $\mathbf{W}^1(G)$ и
   $\mathcal{C}^\infty(\overline{G}
   )$. Из соотношения
       \begin{equation}
\label{rpm__1_}
                 ({rot }+\lambda\,{I})
                 ({rot}- \lambda\,{ I})\mathbf{ u}=-
                 \Delta \,\mathbf{u}
                 +\nabla\, {div}\,\mathbf{ u}-\lambda^2\,\mathbf{ u}
                 \end{equation}
и определения  пространства $V^0(G)$
 видим, что  собственные  функции ротора
  $\mathbf{u}^{\pm}_{\lambda}$,
  %\in                 ,\ne 0
  отвечающие ненулевым собственным значениям
                                 $\pm\lambda $,
   является также собственными  функциями оператора      Лапласа:
                   \begin{equation}
\label{ladi__2_}   -\Delta \mathbf{u}=\lambda^{2}\mathbf{u},\quad
\mathbf{ u}\in V^0(G).
\end{equation}
Нормированные собственные  функции ротора обозначим через
$\mathbf{q}^{\pm}_{j}$.
\[ S\mathbf{q}^{\pm}_{j}=\mathrm{rot}\,\mathbf{q}^{\pm}_{j}=\pm\lambda_j\mathbf{q}^{\pm}_{j}\quad
  \text{при}\,\, \lambda_j\in \Lambda\subset
R, \quad \lambda_j\leq\lambda_{j+1},
\quad\|\mathbf{q}^{\pm}_{j}\|=1.\] Они составляют полный
ортонормированный базис в пространстве $\mathbf{V}^0(G)$.
Спектральное разложение вектор-функции $\mathbf{f}\in V^0(G)$ по
этому базису имеет вид:
\begin{equation}
\label{sp__1_}\mathbf{f}=\sum_{\lambda_j\in
\Lambda}[(\mathbf{f},\mathbf{q}^{+}_{j})\mathbf{q}^{+}_{j}+
(\mathbf{f},\mathbf{q}^{-}_{j})\mathbf{q}^{-}_{j}], \quad
\mathbf{f}\in\mathbf{V}^{0}(G).\end{equation}

В случае шара собственные числа ротора суть корни квадратные из
собственных чисел оператора Лапласа-Дирихле, а собственные функции
ротора вычисляются  явно: \eqref{vrsff   2}%через бесселевы и сферические функции.

Наряду с оператором $S$ %из \cite{giyo},
рассмотрим оператор
%$S+\lambda\,I$:
\begin{equation}
\label{Op__2_} S+\lambda \,I :\mathbf{V}^{0}(G)\rightarrow
\mathbf{V}^{0}(G),\end{equation} который на $\mathbf{W}^{1}(G)$
совпадает с $\mathrm{rot}+\lambda \,I$.

 Оператор $S+\lambda\,I$ является самосопряженным, так как
\begin{equation}
\label{sic__1_} \int_G (\mathrm{rot}+\lambda \,I )\mathbf{u}\cdot
\mathbf{v}\, d\mathbf{x}=\int_{G}\mathbf{u}\cdot
(\mathrm{rot}+\lambda \,I )\mathbf{v} \, d \mathbf{x}
\end{equation} для любых функций  $\mathbf{u}$ и $\mathbf{v}$
из %области определенияоператора $S+\lambda \,I$.
$\mathbf{W}^1(G)$.
 Это доказано в  общем случае в  \cite{giyo}, а в случае шара другим
  способом -  в работе
 автора \cite{saUMJ13}.

 Условие обратимости оператора $S+\lambda\,I$ совпадает с условием:
\begin{equation}
\label{cos__1_} \int_G \mathbf{f}\cdot \mathbf{v}\, dx=0\quad
\forall v\in Ker(S+\lambda\,I). \end{equation}

  Пусть $\mathbf{f}\in V^0(G)$, так как $(S+\lambda\,I)\mathbf{f}\in V^0(G)$
то %этому базису имеет вид:
\begin{equation}
\label{sp__2_} (S+\lambda\,I)\mathbf{f}=\sum_{\lambda_j\in
\Lambda}[(\lambda+\lambda_j)(\mathbf{f},\mathbf{q}^{+}_{j})\mathbf{q}^{+}_{j}+
(\lambda-\lambda_j)(\mathbf{f},\mathbf{q}^{-}_{j})\mathbf{q}^{-}_{j}]
\end{equation}и ряд сходится в $\mathbf{L}_{2}(G)$.
Если $\lambda$ совпадает с одним из собственных значений
$\pm\lambda_{j_0}$, то соответствующее слагаемое в этом ряду
исчезает.

Если элемент $(S+\lambda\,I)^{-1}\mathbf{f}\in V^0(G)$, то
\begin{equation}
\label{sp__3_} (S+\lambda\,I)^{-1}\mathbf{f}=\sum_{\lambda_j\in
\Lambda}[(\lambda+\lambda_j)^{-1}(\mathbf{f},\mathbf{q}^{+}_{j})\mathbf{q}^{+}_{j}+
(\lambda-\lambda_j)^{-1}(\mathbf{f},\mathbf{q}^{-}_{j})\mathbf{q}^{-}_{j}]
\end{equation} и ни одно из слагаемых этого ряда не обращается в
бесконечность. Это означает, что $(\mathbf{f},\mathbf{q}^{-}_{j})=0$
при $\lambda=\lambda_j=\lambda_{j_0}$, то-есть функция $\mathbf{f}$
ортогональна всем собственным функциям
$\mathbf{q}^{-}_{j}(\mathbf{x})$ ротора, отвечающим собственному
значению  $\lambda_{j_0}$.

\begin{theorem} Оператор  $S+\lambda\,I:
\mathbf{V}^{0}(G)\rightarrow \mathbf{V}^{0}(G)$
   однозначно обратим, если $\lambda$ не совпадает ни с одним из
   собственных   значений оператора $S$, и его обратный задается формулой
   \eqref{sp__3_}.

  Если $\lambda=\lambda_{j_0}$, то он обратим тогда и только
  тогда,  когда
  \begin{equation}
\label{urz _1_}\int_G \mathbf{f}\cdot \mathbf{q_j^-}\, dx=0\quad
\text{для}\quad\forall \mathbf{q_j^-}: \lambda_j=\lambda_{j_0}.
\end{equation}
Ядро оператора $S+\lambda_{j_0}\,I$ определяется собственными
функциями $\mathbf{q_j^-}(\mathbf{x})$, собственные значения которых
равны $\lambda_{j_0}$:
\begin{equation} \label{ker__1_}
Ker(S+\lambda_{j_0}\,I)= \sum_{\lambda_j=\lambda_{j_0}}
c_j\,\mathbf{q}^{-}_{j}(\mathbf{x})\quad\text{для}\quad\forall\,\,c_j\in
\mathcal{R}.\end{equation} %где $$--про\cite{bs}извольные постоянные.
\end{theorem}
  Построение собственных функций
ротора в заданной области -- сложная задача. Один случай выделяется
особо.
%%%%p
\section {Построение собственных функций  ротора в шаре}

 Обозначим через $v(\mathbf{x})$ скалярное
произведение векторов $\mathbf{x}$
 и $\mathbf{u}$.
Автор заметил, что внутри шара
 функция $v(\mathbf{x})=\mathbf{x}\cdot\mathbf{u}$ удовлетворяет
уравнению $-\Delta
 v(\mathbf{x})=\lambda^{2}v(\mathbf{x})$,  краевому условию
 $ v|_{S}=0$, и условию  $\quad v(0)=0$ в его центре. Тем самым,
{\it  Любому решению
 $(\lambda,\mathbf{u})$ задачи 1 в шаре $B$ при $\lambda\neq 0 $ соответствует
  решение $(\lambda^{2},\mathbf{x}\cdot \mathbf{u})$ задачи:}
\begin{problem} Найти  собственные значения $\mu $ и собственные
функции $v(x)$  оператора Лапласа $-\Delta $ в шаре $B$ такие, что
      \begin{equation}
\label{ldo__1_}           -\Delta v=\mu\,v \quad \text{в} \quad B,
\quad v|_{S}=0,   \quad v(0)=0.
\end{equation}\end{problem}
Результаты этого параграфа подробно изложены в работе
  \cite{saUMJ13}. Здесь мы
   приведем основные моменты доказательства и
   формулы ее решений.
\subsection{ Функции $\psi_n(z)$.}
\[ \psi_n(z)\equiv\sqrt{\frac{\pi}{2z}}J_{n+\frac12}(z)
 =\sqrt{\frac\pi{2z}}\sum\limits_{p=0}\limits^\infty
 \frac{(-1)^p}{p!\Gamma\bigl(n+1+p+\frac12\bigr)}
 \left(\frac z2\right)^{n+2p+\frac12}.\]
 Как показал Л.~Эйлер (см.\,\cite{vla}, \S 23, с.~356) цилиндрические
 функции $J_{n+\frac12}(z)$ полуцелого порядка
 выражаются через элементарные и
\begin{equation}\label{bes  2}
 \psi_n(z)
 = (-z)^n\left(\frac d{zdz}\right)^n\left(\frac{\sin z}z\right).\\
\end{equation}
 Откуда видно, что
 нули функций $\psi_n(z)$ лежат на действительной оси
 и располагаются на ней симметрично относительно точки $z=0$.

\subsection{Спектральная задача Дирихле для уравнения Лапласа.} В шаре она
  решена в сферической системе
координат $(r,\theta ,\varphi )$ методом разделения переменных
 (см.\cite{vla},  $\S 26$).

 {\it собственные значения
            оператора  $\mathcal{L}$ задачи %в шаре    $B$
            равны $\lambda _{n,m}^{2}=(\rho _{n,m}/ {R})^2$,  где
            $n\ge 0$,   $m\in N$,
         а   числа $\rho _{n,m}>0$ - нули функций $\psi_n(z)$,}

 {\it их действительные собственные функции $v_{\kappa
}^{{}}$ %(соответствующие $\lambda _{n,m}^{2}$)
 имеют вид:
   \begin{equation}
  \label{sfu   1}
    v_{\kappa }(r,\theta ,\varphi )
=c_{\kappa }{\psi }_{n}(\lambda _{n,m}r)Y_{n}^{k}(\theta ,\varphi),
\end{equation}
где $\kappa =(n, m, k)$- мультииндекс, $n\ge 0$, $|k|\le n,$ $m\in
\mathbb{N}$, $c_{\kappa }$-произвольные  постоянные,
 $ P_{n}^{k}(\cos \theta )$ - присоединенные
  функции Лежандра,  $0<r\le R$, $ 0\le \theta \le \pi$,
  $0\le \varphi \le 2\pi$, \quad $Y_{n}^{k}(\theta,\varphi )$-
   сферические функции:}
 \begin{equation}  \label{sffu   2}Y_{n}^{k}(\theta ,\varphi )=
 \left\{\begin{array}{ccc}
  P_{n}^{k}(\cos \theta )\cos(k \varphi), \quad \text{если}
  \quad k=0,1,...,n;\\
P_{n}^{|k|}(\cos \theta )\sin(|k| \varphi), \quad \text{если}\quad
k=-1,...,-n  \end{array}\right.
  \end{equation}
 \subsection{Решение задачи 2.}
  Так как ${{\psi }_{0}}(0)=1$,
 функции $\{v_{\kappa }\}$ при $\kappa =(0,m,0)$
  удовлетворяют   условию $v_{\kappa }^{{}}(0)=0$ задачи 2
    тогда  и только тогда, когда  коэффициенты
    $c_{(0,m,0)}=0$. Откуда следует, что серия $n=0$ в
\eqref{sfu   1} выпадает.

\subsection { Решение задачи 1} %Имеет место\begin{lemma}
В шаре $B$ любому решению  $(\mu ,v)$ задачи 2
 при $\mu > 0 $
соответствуют два и только два решения $(\sqrt{\mu
},{{\mathbf{u}}^{+}})$ и $(-\sqrt{\mu },{{\mathbf{u}}^{-}})$ задачи
1 такие , что $\mathbf{x}\cdot {{\mathbf{u}}^{+}}=\mathbf{x}\cdot
{{\mathbf{u}}^{-}}=v$. %  но ${{\mathbf{u}}^{+}}\ne{{\mathbf{u}}^{-}}$.
 \cite{saUMJ13}.
 Ее обственные значения $\pm\lambda_{n,m}$ - это  корни квадратные из
   собственных чисел задачи 2.

\subsection{Формулы решений задачи} {\it   Ненулевые собственные значения
$\lambda _{n,m}^{\pm }$ задачи 1 равны   $\pm \lambda _{n,m}=\pm
 (\rho_{n,m})/R$,  , где $R$--радиус шара, а числа $\rho
_{n,m}$ -- нули функций $\psi_n(z)$. Собственные функции $u_{\kappa
}^{\pm }$ задачи 1 в сферических координатах вычисляются по
формулам:
\begin{equation}
  \label{vrsff   2}\begin{array}{c}
  u_{\kappa }^{\pm }=c_{\kappa }^{\pm }(\pm\lambda _{n,m}
r)^{-1}{\psi }_{n}
  (\pm\lambda _{n,m}r)Y_{n}^{k}(\theta ,\varphi  )\,\mathbf{i}_r+\\
c_{\kappa }^{\pm }{{(\pm\lambda _{n,m}^{\pm }r)}^{-1}}
   Re[\Phi _{n}(\pm\lambda _{n,m}r)](Re HY_{n}^{k}\,
  \mathbf{i}_\varphi+
  Im HY_{n}^{k}\,\mathbf{i}_\theta)+\\
  c_{\kappa }^{\pm }{{(\pm\lambda _{n,m}r)}^{-1}}
  Im[\Phi _{n}(\pm\lambda _{n,m}r)](-Im HY_{n}^{k}\,
  \mathbf{i}_\varphi+
  Re HY_{n}^{k}\,\mathbf{i}_\theta).\end{array}\end{equation}
                где числа %$i$ - мнимая единица,
  $c_{\kappa }^{\pm }\in \mathbb{R}$,
$m{{,}^{{}}}n\in \mathbb{N}$,    $|k|\le n$,$\kappa=(n,m,k)$,
 $Y_{n}^{k}(\theta ,\varphi  )$--сферические функции, $\mathbf{i}_r, \mathbf{i}_\theta,
\mathbf{i}_\varphi$-репер,}
     \[\Phi _{n}^{{}}(\pm\lambda _{n,m}r)=
  \overset{{}}{\mathop{{}}}\,\int\limits_{0}^{r}{}\,{{e}^{\pm\,i\lambda _{n,m}
  (r-t)}}^{{}}
  {{\psi }_{n}}{{(\pm\lambda _{n,m}t)}}{{t}^{-1}}dt,
  \quad Im \Phi _{n}(\pm\rho_{n,m})=0,\]%\end{equation}
                                        \begin{equation}
  \label{oph   1}
  \text{H}Y_{n}^{k}(\theta ,\varphi )=
  {{\left( {{\sin }^{-1}}{{\theta }^{_{{}}}}{{\partial }_{\varphi }}+
  i{{\partial }_{\theta }} \right)}^{{}}}Y_{n}^{k}(\theta ,\varphi ).   \end{equation}
   Эти формулы используются при рассчетах поля скоростей $u_{\kappa }^{\pm
}(x)$ вихревого потока  при заданном $\kappa$.

Г.Г.Исламов (Удмурдский ГУ, Ижевск), используя программы Volfram
Mathematica рассчитал эти поля при минимальном собственном значении
и траектории их линий тока. %\footnote
({\small см. его доклад  в %материалах конференции
http://www.wolfram.com/events/technology-conference-ru/2016/
resources.html})

Траектория отдельной точки похожа на нить, которая наматывается на
тороидальную катушку, каждая на свою.

 В работе \cite{cdtgt} также
определены собственные функции ротора при минимальном собственном
значении
 и рассчитаны траектории их линий тока.

  Идея  сведения краевой задачи
$\mathbf{n}\cdot\mathbf{u}|_{S}=g$
 для системы $\mathbf{rot}\mathbf{u}+\lambda\mathbf{u}=\mathbf{f}$
 в шаре при $\lambda\neq 0$  к задаче Дирихле для
уравнения Гельмгольца возникла давно \cite{sa71, sa72}. {\footnote
{\small В 1970 А.А.Фурсенко,  студент НГУ,   в дипломной работе
таким способом решил эту задачу   в классах Гельдера.
 Мы   выписали формулы  решений  задачи, но  не
 опубликовали их. Я опубликовал  формулы \eqref{vrsff   2}
   в 2000 году \cite{sa2000, sa01}, когда узнал о приложениях и о работе
 \cite{chake}.}}

\subsection{Сходимость ряда
 Фурье по собственным функциям
 ротора в норме пространства Соболева $\mathbf{H}^s(B)$, $s\geq 1$} Положим
  \[\mathbf{V}^s_\mathcal{R}(B)=
\{\mathbf{f}\in \mathbf{V}^0\cap\mathbf{H}^s(B):
\mathbf{n}\cdot\mathbf{f}|_S=0,...,
 \mathbf{n}\cdot\ \text{rot}^{s-1}\mathbf{f}|_S=0,\,\,
 \|\mathbf{f}\|_{\mathbf{V}_\mathcal{R}^s}=\|\mathbf{f}\|_{\mathbf{H}^s}\}.\]
 \begin{theorem}
 Для того, чтобы $\mathbf{f}\in \mathbf{V}^0(B)$
 разлагалась в ряд Фурье
 \begin{equation} \label{rof 1}
\mathbf{f}(\mathbf{x})=\sum_{\kappa, n>0}
((\mathbf{f},\mathbf{q}_{\kappa}^+)\mathbf{q}_{\kappa}^+(\mathbf{x})
+(\mathbf{f},\mathbf{q}_{\kappa}^-)\mathbf{q}_{\kappa}^-(\mathbf{x})),
\quad \|\mathbf{q}_{\kappa}^{\pm}\| =1,
\end{equation}
 по  собственным вектор-функциям $\mathbf{q}_{\kappa}^{\pm}(\mathbf{x})$
 ротора в шаре,
 сходящийся в норме
 пространства Соболева $\mathbf{H}^s(B)$, необходимо и достаточно,
  чтобы $\mathbf{f}$ принадлежала $\mathbf{V}^s_\mathcal{R}(B)$.

 Если $\mathbf{f}\in \mathbf{V}^s_\mathcal{R}(B)$,
то сходится ряд
\begin{equation} \label{rof 2}
\sum_{\kappa, n>0}{\lambda}_{\kappa}^{2s}\,
(|(\mathbf{f},\mathbf{q}_{\kappa}^+)|^2
+|(\mathbf{f},\mathbf{q}_{\kappa}^-|^2)),\quad
{\lambda}_{\kappa}=({\rho}_{n,m})/R
\end{equation} и существует такая положительная постоянная $C>0$, не
зависящая от $\mathbf{f}$, что
\begin{equation} \label{orf 3}
\sum_{\kappa, n>0} {\lambda}_{\kappa}^{2s}\,
(|(\mathbf{f},\mathbf{q}_{\kappa}^+)|^2
+|(\mathbf{f},\mathbf{q}_{\kappa}^-|^2))\leq
C\|\mathbf{f}\|^2_{\mathbf{H}^s(B)}.
\end{equation}

 Если $s\geq 2$, то любая вектор-функция  $\mathbf{f}$ из
  $\mathbf{V}^s_\mathcal{R}(B)$
разлагается в в ряд Фурье, сходящийся в пространстве $\mathbf{C}^{s-2}(\overline{B})$.%
\end{theorem}

Действительно, граница шара $S\in \mathcal{C}^\infty$ и собственные
 функции  ${q}_{\kappa}^{\pm}(\mathbf{x})$, $\text{rot}\,{q}_{\kappa}^{\pm}(\mathbf{x})=
 \pm{\lambda}_{\kappa}{q}_{\kappa}^{\pm}(\mathbf{x})$, первой краевой задачи для
 оператора ротор в шаре принадлежат классу
 $\mathcal{C}^{\infty}$  в  $\overline{B}$ . Значит, они
 и их конечные линейные комбинации %\newline
 $\sum_{\kappa}
(c^+_{\kappa}\,{q}^+_{\kappa}(\mathbf{x})+c^-_{\kappa}\,{q}^-_{\kappa}(\mathbf{x}))$
принадлежат любому из пространств  $\mathbf{V}^l_\mathcal{R}(B)$ при
$l>0$ и
$\gamma\mathbf{n}\cdot\mathbf{q}^{\pm}_{\kappa}(\mathbf{x})=0$,\newline
 $\gamma\mathbf{n}\cdot \text{rot}\,\mathbf{q}^{\pm}_{\kappa}(\mathbf{x})=0$,
и так далее.

 Для следа нормальной компоненты  вектор-функции
 $\mathbf{f}\in  \mathbf{H}^l(B)$ на $S$ и
 ее производных
$\gamma\partial^\alpha \mathbf{f}$ при $|\alpha|<l$ имеются оценки
(см. \cite{mi}, $\S\ 5.1 $ главы 3):
\begin{equation}
  \label{onaf  2}
  \|\gamma\mathbf{n}\cdot\partial^{\alpha} \mathbf{f}\|_{L_2(S)}\leq
  \|\gamma\partial^{\alpha} \mathbf{f}\|_{\mathbf{L}_2(S)}\leq c \| \mathbf{f}|\|_{\mathbf{H}^{|\alpha|+1}(B)}
  \leq c \|\mathbf{f}|\|_{\mathbf{H}^{l}(B)}.
  \end{equation}

 Обозначим через $\mathbf{S}_N(\mathbf{x})$ частичную сумму ряда
 \eqref{rof 1},   $\mathbf{S}_N(\mathbf{x})\in \mathbf{V}^s_\mathcal{R}(B)$ при всех
$N\geq 1$ и $s\geq 1$. Рассмотрим разность
$\mathbf{f}(\mathbf{x})-\mathbf{S}_N(\mathbf{x})$ и воспользуемся
оценкой %\eqref{onaf 2}
 следа   на $S$ ее нормальной компоненты:
\begin{equation}
  \label{oaf  0}
  \|\gamma\mathbf{n}\cdot (\mathbf{f}-\mathbf{S}_N)\|_{L_2(S)}\leq
   \|\gamma (\mathbf{f}-\mathbf{S}_N)\|_{\mathbf{L}_2(S)}\leq
  c \| (\mathbf{f}-\mathbf{S}_N)|\|_{{\mathbf{H}}^{1}(B)}
    \end{equation}
 Если функция $\mathbf{f}\in  \mathbf{V}^0(B)$ и
ряд Фурье \eqref{rof 1}сходится в норме $\mathbf{H}^1(B)$, то
$\|\mathbf{f}(\mathbf{x})-\mathbf{S}_N(\mathbf{x})\|_{{\mathbf{H}}^{1}(B)}\rightarrow
0$ при $N\rightarrow \infty$. Так как
$\gamma\mathbf{n}\cdot\mathbf{S}_N=0$ при любых $N$, то $ \|
\gamma\mathbf{n}\cdot\mathbf{f}\|_{L_2(S)}=0$ и, значит, $
\gamma\mathbf{n}\cdot\mathbf{f}=0$ и $\mathbf{f}\in
\mathbf{V}_{\mathcal{R}}^1(B)$.

 Если
 функция $\mathbf{f}\in  \mathbf{V}^0(B)$ и
ряд Фурье \eqref{rof 1} сходится в норме $\mathbf{H}^2(B)$, то
воспользуемся оценкой следа    нормальной компоненты ротора:
%\begin{equation}  \label{oaf  0}
 \[ \|\gamma\mathbf{n}\cdot \text{rot}\,(\mathbf{f}-\mathbf{S}_N)\|_{L_2(S)}\leq
   \|\gamma\text{rot}\, (\mathbf{f}-\mathbf{S}_N)\|_{\mathbf{L}_2(S)}\leq
  c \| (\mathbf{f}-\mathbf{S}_N)|\|_{{\mathbf{H}}^{2}(B)}.\]
    %\end{equation}
Так как $\gamma\mathbf{n}\cdot \text{rot}\,\mathbf{S}_N=0$ при любых
$N$, то аналогично     предыдущему
 $\gamma\mathbf{n}\cdot\mathbf{f}=0$  и
  $\gamma\mathbf{n}\cdot \text{rot} \mathbf{f}=0$.
 Значит, $\mathbf{f}\in
\mathbf{V}_{\mathcal{R}}^2(B)$.\newline
 И так далее, если
функция $\mathbf{f}\in \mathbf{V}^0(B)$  и ряд Фурье \eqref{rof 1}
сходится в норме $\mathbf{H}^s(B)$, то $\mathbf{f}\in
\mathbf{V}^s_\mathcal{R}(B)$, где $s>2$ . Необходимость доказана.

 Пусть $\mathbf{f}\in \mathbf{V}^s_\mathcal{R}(B)$,
где $s>0$ . Установим справедливость неравенства \eqref{orf 3}. Так
как $\text{rot}\, \mathbf{q}^{\pm}_{\kappa}={\pm}\lambda_\kappa
\mathbf{q}^{\pm}_{\kappa}$, согласно формуле Грина %\eqref{rot 1}
имеем
\begin{equation}
\label{rot   2}  \int\limits_{B}\,
\text{rot}\,{\mathbf{f}}\cdot\mathbf{q}^{\pm}_{\kappa}
\,d{\mathbf{x}}= {\pm}\lambda_\kappa\int\limits_{B}{\mathbf{f}}
\cdot\mathbf{q}^{\pm}_{\kappa}\, d{\mathbf{x}}+
\int\limits_{S}[\mathbf{f},\mathbf{q}^{\pm}_{\kappa}]\cdot{\mathbf{n}}\,
dS.
\end{equation} Сокращенно эту формулу запишем так
\begin{equation}
\label{rot   3} (\text{rot}\,{\mathbf{f}},\mathbf{q}^{\pm}_{\kappa})
= {\pm}\lambda_\kappa({\mathbf{f}} ,\mathbf{q}^{\pm}_{\kappa}) +
{c}^{\pm}_{\kappa}(\mathbf{f}),\quad {c}^{\pm}_{\kappa}(\mathbf{f})=
\int\limits_{S}[\mathbf{f},\mathbf{q}^{\pm}_{\kappa}]\cdot{\mathbf{n}}\,
dS. \end{equation}
%\begin{equation}%\end{equation}
%\label{of   1} \end{equation}%\begin{equation}\label{nof   1}
 ${c}^{\pm}_{\kappa}(\mathbf{f})$-
ограниченный функционал над $\mathbf{H}^1(B)$, так как
\[|{c}^{\pm}_{\kappa}(\mathbf{f})|\leq \int\limits_{S}|\mathbf{f}|
|\mathbf{q}^{\pm}_{\kappa}|\, dS\leq
\|\mathbf{f}\|_{\mathbf{L}_2(S)}\,\|\mathbf{q}^{\pm}_{\kappa}\|
_{\mathbf{L}_2(S)}\leq\|\mathbf{f}\|_{1}\,\|\mathbf{q}^{\pm}_{\kappa}\|
=\|\mathbf{f}\|_{1}.\]

Отметим, что ${c}^{\pm}_{\kappa}(\mathbf{f})=0$, если
$\mathbf{f}|_S=0$ или если $\mathbf{f}=\mathbf{S}_N$ при $N<\infty$.

 Обозначим через $\beta_{\kappa}^{\pm}$
коэффициенты Фурье функции $\text{rot}^s \mathbf{f}$. Согласно
формуле \eqref{rot 3}
\begin{equation}
\label{rot   s}\beta_{\kappa}^{\pm}=(\text{rot}^s\, \mathbf{f},
\mathbf{q}^{\pm}_{\kappa})={\pm}\lambda_\kappa(\text{rot}^{s-1}\,{\mathbf{f}}
,\mathbf{q}^{\pm}_{\kappa}) +
{c}^{\pm}_{\kappa}(\text{rot}^{s-1}\,\mathbf{f})=...\end{equation}
\[({{\pm}\lambda_\kappa})^{s}(\,{\mathbf{f}}
,\mathbf{q}^{\pm}_{\kappa})+
({\pm}\lambda_\kappa)^{s-1}\,{c}^{\pm}_{\kappa}(\mathbf{f})+
({\pm}\lambda_\kappa)^{s-2}\,{c}^{\pm}_{\kappa}(\text{rot}\,\mathbf{f})+...
+{c}^{\pm}_{\kappa}(\text{rot}^{s-1}\,\mathbf{f}).\]

 Поскольку $\text{rot}^s \mathbf{f}\in
\mathbf{L}_2(B)$, то
 \begin{equation}
  \label{rodf  1}\Sigma_\kappa [(\beta_{\kappa}^+)^2+
(\beta_{\kappa}^-)^2]= \|\text{rot}^s \mathbf{f}\|^2.\end{equation}
Для финитных вектор-функций из %$\mathbf{f}\in
$\mathbf{V}^s_\mathcal{R}(B)$ имеем
\begin{equation} \label{orf 4}
\sum_{\kappa, n>0} {\lambda}_{\kappa}^{2s}\,
(|(\mathbf{f},\mathbf{q}_{\kappa}^+)|^2
+|(\mathbf{f},\mathbf{q}_{\kappa}^-|^2))= \|\text{rot}^s
\mathbf{f}\|^2\leq C\|\mathbf{f}\|^2_{\mathbf{H}^s(B)}.
\end{equation}
Но финитные вектор-функции  $\mathbf{f}\in \mathcal{D}(B)$ из
$\mathbf{V}^s_\mathcal{R}(B)$ плотны в
$\mathbf{V}^s_\mathcal{R}(B)$. Неравенство \eqref{orf 3} доказано.

Вернемся к   частичной сумме $\mathbf{S}_l(\mathbf{x})$ ряда
\eqref{rof 1}. Как мы уже отмечали $\mathbf{S}_l(\mathbf{x})\in
\mathbf{V}^s_\mathcal{R}(B)$ при всех $l>0$. В частности,
$\text{div} \mathbf{S}_l(\mathbf{x})=0$ и \newline
$\gamma\mathbf{n}\cdot\mathbf{S}_l(\mathbf{x})=0$.
 Поэтому оценка
\eqref{apro s} при $s=0$ принимает вид
\begin{equation} \label{apr 1} C_1\|\mathbf{S}_l\|_{1}
\leq\|\mathrm{rot}\,\mathbf{S}_l\|+ \|\mathbf{S}_l\|.\end{equation}

Легко видеть, что $\|\mathbf{S}_l\|^2\leq
c\|\mathrm{rot}\mathbf{S}_l\|^2$, где
$c=max_{m,n}\lambda_{m,n}^{-2}$.
 Поэтому
\begin{equation} \label{apr 2} \|\mathbf{S}_l\|^2_{1}
\leq a_1\|\mathrm{rot}\,\mathbf{S}_l\|^2.\end{equation}
%\|\mathbf{S}_l\|^2). Следовательно,
По индукции при $s>1$
\begin{equation} \label{apr s} \|\mathbf{S}_l\|^2_{s}
\leq a_s\|\mathrm{rot}^s\,\mathbf{S}_l\|^2.\end{equation}
%\|\mathrm{rot}^{s-1}\,\mathbf{S}_l\|^2+...+ \|\mathbf{S}_l\|^2).
Пусть $\mathbf{f}\in \mathbf{V}^s_\mathcal{R}(B)$, где $s>0$.
Согласно неравенству \eqref{orf 3}, ряды в его левой части
 сходятся и если $l>m\geq 1$, то
\[\|\mathbf{S}_l-\mathbf{S}_m\|^2_{s}\leq a_s\|rot^s
(\mathbf{S}_l-\mathbf{S}_m)\|^2=\]
\[a_s\sum_{m+1}^l\lambda_{\kappa}^{2s}(|(\mathbf{f},\mathbf{q}_{\kappa}^+)|^2
+|(\mathbf{f},\mathbf{q}_{\kappa}^-|^2))\rightarrow 0\]
 при
$l,m\rightarrow\infty$. Это означает, что ряд \eqref{rof 1} сходится
к $\mathbf{f}$ в $\mathbf{H}^s(B)$.

 При $s\geq 2$ в трехмерном шаре $B$ имеется
вложение пространств
$\mathbf{H}^s(B)\subset\mathbf{C}^{s-2}(\overline{B})$ и оценка:
\begin{equation}
  \label{cw  1}\|\mathbf{f}\|_{\mathbf{C}^{s-2}(\overline{B})}\leq C_s\|
\mathbf{f}\|_{\mathbf{H}^s (B)}\end{equation} для любой функции
$\mathbf{f}\in \mathbf{H}^s (B)$, в которой постоянная  $C_s>0$ не
зависит от $\mathbf{f}$ (см., например, Теорему 3 $\S\,6.2$ в
\cite{mi})). В частности,
\begin{equation}
  \label{cs  l}\|\mathbf{S}_l-\mathbf{S}_m\|_{\mathbf{C}^{s-2}(\overline{B})}
  \leq C_s\|\mathbf{S}_l-\mathbf{S}_m\|_{\mathbf{H}^s (B)}.\end{equation} Если
  $\|\mathbf{S}_l-\mathbf{S}_m\|_{\mathbf{H}^s (B)}\rightarrow 0$ при
$l,m\rightarrow\infty$, то
$\|\mathbf{S}_l-\mathbf{S}_m\|_{\mathbf{C}^{s-2}(\overline{B})}\rightarrow
0$. Это означает, что ряд \eqref{rof 1} сходится к $\mathbf{f}$ в
$\mathbf{C}^{s-2}(\overline{B})$.
Теорема доказана.%15.11.2013

{\bf Следствие.} {\it  Любая соленоидальная вектор-функция
$\mathbf{f}$ из $\mathbf{C}^{\infty}_0({B})$
 разлагается в
 ряд Фурье \eqref{rof 1},
сходящийся в пространстве $\mathbf{C}^{\infty}(\overline{B})$.}
\subsection{Скалярное произведение функций $\mathbf{f}$ и
$\mathbf{g}$ из $\mathcal{B}$ в базисе из собственных функций
ротора} Оно имеет вид:
\begin{equation} \label{spb 1}( \mathbf{f}, \mathbf{g})=
\sum_{\kappa, n>0} %{\lambda}_{\kappa}^{2s}\,
[(\mathbf{f},\mathbf{q}_{\kappa}^+)(\mathbf{g},\mathbf{q}_{\kappa}^+)
+(\mathbf{f},\mathbf{q}_{\kappa}^-)(\mathbf{g},\mathbf{q}_{\kappa}^-)].
\end{equation}
Если $\mathbf{f}$ и $\mathbf{g}$ принадлежат
$\mathbf{V}^1_\mathcal{R}(B)$, то равенства
\begin{equation} \label{sor }(\text{rot}\, \mathbf{f}, \mathbf{g})=
(\mathbf{f},\text{rot}\mathbf{g})= \sum_{\kappa, n>0}
{\lambda}_{\kappa}[(\mathbf{f},\mathbf{q}_{\kappa}^+)(\mathbf{g},\mathbf{q}_{\kappa}^+)
-(\mathbf{f},\mathbf{q}_{\kappa}^-)(\mathbf{g},\mathbf{q}_{\kappa}^-)
]
\end{equation} показывают, что оператор $\mathbf{rot}$ является
 самосопряженным в $\mathcal{B}$.

    %%%%
\section{ Градиент дивергенции
 в ограниченной области}

\subsection{Краевая задача:}
в ограниченной области $G$ с гладкой границей $\Gamma$ заданы
векторная и скалярная функции $\mathbf{f}$ и ${g}$, найти
вектор-функцию $\mathbf{u}$, такую что
\begin{equation}
\label{grd__2_}
 \nabla\, \text{div}\mathbf{u}+\lambda\, \mathbf{u}=\mathbf{f}
 \quad  \text{в}\quad G,
\quad \mathbf{n}\cdot\mathbf{u}|_\Gamma=g,
\end{equation}
  При  $\lambda\ne 0$ она  является обобщенно
   эллиптической (см. $\S 2$). Действительно,
оператор $\nabla\, \text{div}+\lambda\text{I}$ второго порядка
принадлежит классу (RNS1)   (см. п.2.1 ), так как
 а)  ранг его символической матрицы
$\nabla\, \text{div}(i\xi)$ постоянный и равен единице, б)оператор
$\nabla\,\text{div}$  имеет  левый
 аннулятор $\text{rot}$: $ \text{rot}\,\nabla\, \text{div} \mathbf{u}=0$,
в) оператор
 $\nabla\,\text{div}$// $\lambda\text{rot}$ эллиптичен. % (см.  $1^0$).
 Поэтому расширенная система:
  \begin{equation}
\label{rodi__2_} \nabla\, \text{div}\,\mathbf{u}+\lambda
\mathbf{u}=\mathbf{f},\quad \lambda \text{rot}\,\mathbf{
u}=\text{rot}\,\mathbf{f}\end{equation}является эллиптической
системой при выборе порядков: $s_k=0$  при $k=1,-,3$ и $s_k=-1$ при
$k=4,-,6$; $t_j=2$ при $ j=1,-,3$.

Далее,  система \eqref{rodi__2_} с краевым условием
$\gamma\mathbf{n}\cdot\mathbf{u}=g$
 эллиптична по определению В.А.Солонникова, так как

1) Расширенная  переопределенная система \eqref{rodi__2_}
эллиптична,

2) граничный оператор $\gamma\mathbf{n}\cdot\mathbf{u}$ "накрывает"
\,оператор системы \eqref{rodi__2_}.

А именно, $1^0)$ однородная система линейных алгебраических
уравнений:
\begin{equation}
\label{cdx  2} \lambda\,\text{rot}(i\xi )\mathbf{w}=0, \quad
(\nabla\,\text{div})(i\xi )\mathbf{w}=0, \quad \forall \xi\neq 0
\end{equation}
 c параметром $\xi \in T'( G)$ имеет только тривиальное
 решение $ \mathbf{w}=0$;

  $2^0)$   однородная система линейных
дифференциальных уравнений:
\begin{equation}\label{cdz  3}
\lambda\text{rot}(i\tau+\mathbf{n} d/dz ) \mathbf{v}=0,\quad
(\nabla\text{div})(i\tau+\mathbf{n} d/dz )\mathbf{v}=0, \quad
\forall \tau \neq 0,\end{equation} на полуоси $z\geq 0$ с краевым
условием: $\mathbf{n}\cdot \mathbf{v}|_{ z=0}=0$ и условием
убывания: $\mathbf{v}(y,\tau; z )\rightarrow 0$ при $z\rightarrow +
\infty$, имеет только тривиальное решение. Здесь $\tau$ и
$\mathbf{n}$ касательный и нормалльный векторы к
 $\Gamma$ в точке $y\in \Gamma$ и $|\mathbf{n}|=1$.

Доказательство этих утверждений {\footnote {\small Ввиду ограничения
объема статьи я вынужден опустить  это и другие подобные
рассуждения. Надеюсь, читатель без труда восстановит их. }} почти
такое же, как в п.2.1.

Итак, %система \eqref{rodi__2_} с краевым условием из \eqref{grd__2_}
 %эллиптична. По определению
 краевая задача  \eqref{grd__2_} при $\lambda\neq 0$
  является обобщенно эллиптической.

 %%%5.2
\subsection{Оператор задачи \eqref{grd__2_} в пространствах Соболева}
Пусть $\mathbf{u}$ принадлежит пространству ${\mathbf{H}^{s+2}}(G)$,
то-есть каждая компонента $u_j\in H^{s+2}(G)$. Тогда
$\nabla\text{div} \mathbf{u}$ принадлежит $\mathbf{H}^{s}(G)$, и
$\text{rot}^{2}\mathbf{u}$ принадлежит $\mathbf{H}^{s}(G)$. Поэтому
вектор-функция $\mathbf{f}:=\nabla\text{div} \mathbf{u}+ \lambda
\mathbf{u}$ принадлежит пространству
\begin{equation}\label{pr  2} {\bf{F}^{s}}(G)=\{\mathbf{f}\in
{\mathbf{H}^{s}(G)}: \text{rot}^2\,\mathbf{f}\in
 \mathbf{H}^{s}(G)\},\end{equation} которое снабдим  нормой
 \[\|\mathbf{v}\|_{\mathbf{F}^{s}}=(
 \|\mathbf{v}\|^2_{s}+
 \|\text{rot}^2\mathbf{v}\|^2_{s})^{1/2}.\] Функция
$g:=\gamma({\mathbf{n}}\cdot\mathbf{u})$
%\equiv\mathbf{n}\cdot\mathbf{u}|_{\Gamma}$
принадлежит пространству
$H^{s+3/2}(\Gamma)$. Следовательно, при $\lambda\neq 0$ задаче
соответствует ограниченный оператор
\begin{equation}\label{op  2} \mathbb{B}\mathbf{u}\equiv \begin{matrix}
\nabla\text{div}\,\mathbf{u}+\lambda\,\mathbf{u} \\
\gamma({\mathbf{n}}\cdot\mathbf{u})\end{matrix}:
\bf{H}^{s+2}(G)\rightarrow
\begin{matrix}\bf{F}^{s}(G)\\ H^{s+3/2}(\Gamma)\end{matrix}.
\end{equation}
Согласно Теоремем 1.1 в работе Солонникова \cite{so71}, о
переопределенных эллиптических краевых задачах  в ограниченной
области $G$   с гладкой границей $\Gamma\in \mathcal{C}^{s+2} $,
обобщенно эллиптический оператор \eqref{op 2} имеет левый
регуляризатор, то-есть
 ограниченный оператор  $\mathbb{B}^L$ такой, что
$\mathbb{B}^L\mathbb{B}=\mathbb{I}+\mathbb{T}$, где $\mathbb{I}$ -
единичный, а $\mathbb{T}$ - вполне непрерывный операторы, и
существует постоянная $C_s >0$ такая, что выполняется априорная
оценка:
\begin{equation} \label{apgro s} C_s\|\mathbf{u}\|_{s+2}
\leq|\lambda|\,\|\mathrm{rot}^2\,\mathbf{u}\|_{s}+
\|\nabla\mathrm{div}\,\mathbf{u}\|_{s}+
|\gamma({\mathbf{n}}\cdot\mathbf{u})|_{s+3/2}+
\|\mathbf{u}\|_{s}.\end{equation} {\small А также аналогичная оценка
в  пространствах Соболева    $W^{s+2}_p(G)$, $ p>1$.} Значит, имеет
место
\begin{theorem}
Оператор $\mathbb{B}$ в пространствах (\ref{op  2}) имеет левый
регуля-ризатор. Его ядро $\mathcal{M}$ конечномерно и выполняется
 оценка \eqref{apgro s}.
 \end{theorem}

Из этой теоремы и  оценки следует, что при $\lambda \ne 0$

 a){\it число линейно независимых решений однородной
задачи \eqref{grd__2_} конечно,}

 b){\it любое ее обобщенное решение %задачи
 бесконечно дифференцируемо вплоть до
границы, если граница области бесконечно дифференцируема.}
\subsection{Оператор $\nabla\mathrm{div}+\lambda I$ в подпространствах }
На подпространстве $\mathcal{B}$ в $\mathbf{L}_2(G)$, ортогональном
подпространству
  $\mathcal{A}$, %=$\{\mathbf{u}=\nabla\, h:\,\, h\inH^1(G)\}$
 оператор $\nabla\, \text{div}\mathbf{u} +\lambda \mathbf{u} $
является  оператором умножения: $\lambda \mathbf{u} $.

 Пространство
$\mathcal{A}_{\gamma}=\{\mathbf{u}=\nabla\, h:\,\, h\in H^1(G),\quad
(\mathbf{n}\cdot\mathbf{u})|_{\Gamma}=0 \}$ плотно в $\mathcal{A}$,
так как функции из $\mathcal{C}_0^\infty \cap \mathcal{A}_{\gamma}$
плотны в $\mathbf{L}_2(G)$. Пространство
  \begin{equation}
\label{sa__1_}\mathcal{A}^2(G)=\{ \mathbf{v}\in
\mathcal{A}_{\gamma}: \nabla\mathrm{div} \mathbf{v}\in
\mathcal{A}_{\gamma}\}.\end{equation}  плотно в
$\mathcal{A}_{\gamma}$ и
 содержится в $\mathbf{H}^2(G)$ в силу оценки \eqref{apgro s}.

Введем оператор $\mathcal{N}_d:\mathcal{A}_{\gamma}
 \rightarrow \mathcal{A}_{\gamma}$ с  областью
определения   $\mathcal{A}^2(G)$, %  : Оператор $\mathcal{N}_d$
который совпадает с $\nabla\mathrm{div}\,\mathbf{v}$ при
$\mathbf{v}\in \mathcal{A}^2(G)$.

Оператор $\mathcal{N}_d+\lambda I:\mathcal{A}_{\gamma}\rightarrow
\mathcal{A}_{\gamma}$ является самосопряженным (эрмитовым).
 Действительно, согласно формуле Гаусса-Остроградсного
\begin{equation}
\label{gri__2_} \int_G(\nabla\mathrm{div}\mathbf{u}+\lambda
\mathbf{u})\cdot\,\mathbf{v} dx=
\int_G\mathbf{u}\cdot\,(\nabla\mathrm{div}\mathbf{v}+\lambda
\mathbf{v})dx+\end{equation}%\mathbf{u},
\[\int_{\Gamma}[(\mathbf{n}\cdot\mathbf{v})\mathrm{div}\mathbf{u}+
(\mathbf{n}\cdot\mathbf{u})\mathrm{div}\mathbf{v}]|_\Gamma \, d S.\]

Если вектор-функции $\mathbf{u}$ и $\mathbf{v}$ принадлежат
$\mathcal{A}^2(G)$, то граничные интегралы пропадают, остальные
интегралы сходятся. Следовательно,
\begin{equation}
\label{gri__2_} ((\nabla\mathrm{div}\mathbf{u}+\lambda
\mathbf{u}),\mathbf{v})=
(\mathbf{u},(\nabla\mathrm{div}\mathbf{v}+\lambda \mathbf{v})) \quad
\text{в}\quad \mathbf{L}_2(G).\end{equation}

%%%%%

{\it Область определения $\mathcal{A}^2(G)$ оператора
$\mathcal{N}_d$ содержится в $\mathbf{H}^2(G)$ и плотна в
$\mathcal{A}_{\gamma}$, а  область его значений совпадает с
$\mathcal{A}_{\gamma}$.}
%\newline

Так как $\mathcal{A}_\gamma$ ортогонально $Ker\mathcal{N}_d$,
оператор $\mathcal{N}_d$ имеет единственный обратный
$\mathcal{N}_d^{-1}$ определенный на $\mathcal{A}_\gamma$. Оператор
%$\mathcal{N}_d$ замкнут и
 %имеет компактный обратный
 $\mathcal{N}_d^{-1}:\mathcal{A}_{\gamma}\rightarrow
 \mathcal{A}^2(G)$ имеет точечный спектр, который
 не содержит точек накопления кроме нуля
Следовательно,  спектр самосопряженного  оператора $\mathcal{N}_d$
точечный  и действительный, а система его собственных вектор-функций
ортогональна и полна в пространстве $\mathcal{A}_{\gamma}$. Каждому
собственному значению соответствует конечное число  собственных
вектор-функций.

Пусть $\mathbf{f}\in \mathcal{A}_{\gamma}(G)$, так как
 $(\mathcal{N}_d+\lambda\,I)\mathbf{f}\in \mathcal{A}_{\gamma}(G)$
то
\begin{equation}
\label{sp__2_} (\mathcal{N}_d+\lambda\,I)\mathbf{f}=\sum_{\mu_j\in
{M}}[(\lambda+\mu_j)(\mathbf{f},\mathbf{q}_{j})\mathbf{q}_{j} ]
\end{equation}и ряд сходится в $\mathbf{L}_{2}(G)$.
Если $\lambda+\mu_{j_0}=0$,
 то соответствующее слагаемое в этом ряду исчезает.

Если элемент $(\mathcal{N}_d+\lambda\,I)^{-1}\mathbf{f}\in
\mathcal{A}_{\gamma}(G)$, то
\begin{equation}
\label{sp__4_}
(\mathcal{N}_d+\lambda\,I)^{-1}\mathbf{f}=\sum_{\mu_j\in
M}[(\lambda+\mu_j)^{-1}(\mathbf{f},\mathbf{q}_{j})\mathbf{q}_{j} ]
\end{equation} и ни одно из слагаемых этого ряда не обращается в
бесконечность. Это означает, что $(\mathbf{f},\mathbf{q}_{j})=0$ при
$-\lambda=\mu_j=\mu_{j_0}$, то-есть функция $\mathbf{f}$
ортогональна всем собственным функциям $\mathbf{q}_{j}(\mathbf{x})$
градиента дивергенции, отвечающим собственному значению $\mu_{j_0}$.
 Итак, имеет место
\begin{theorem}
a). Оператор $\mathcal{N}_d:\mathcal{A}_{\gamma}\rightarrow
\mathcal{A}_{\gamma}$   является
 самосопряженным.
  Его спектр  $\sigma(\mathcal{N}_d )$
 точечный  и действительный.
  Семейство собственных функций $\mathbf{q}_{j}(x)$ оператора
  $\mathcal{N}_d$ образует полный ортонормированный базис в  пространстве
   $\mathcal{A}_{\gamma}$;
    разложение   $\mathbf{a}(x)\in{\mathcal {{A}} _{\gamma}} (G)$ имеет вид
\begin{equation}
\label{spr__4_}\mathbf{a}(x)=\sum_{\mu_j\in
M}(\mathbf{a},\mathbf{q}_{j}) \mathbf{q}_{j}(x),\quad
\|\mathbf{q}_{j}\|=1.\end{equation}% \label{sp__4_}
  b). Если $-\lambda$ не совпадает ни с одним из
   собственных   значений оператора $\mathcal{N}_d$, то
   оператор $\mathcal{N}_d+\lambda\,I:
    \mathcal{A}_{\gamma}\rightarrow\mathcal{A}_{\gamma}$ однозначно обратим,
    и его обратный   задается формулой   \eqref{sp__4_}.
Если $-\lambda=\mu_{j_0}$, то он обратим тогда и только
  тогда,когда
  \begin{equation}
\label{urz _1_}\int_G \mathbf{f}\cdot \mathbf{q_j}\, dx=0\quad
\text{для}\,\,\forall \mathbf{q_j}: %\quad\text{у которых}\,\,
\mu_j=\mu_{j_0}.
\end{equation}
Ядро оператора $\mathcal{N}_d-\mu_{j_0}\,I$ определяется
собственными функциями $\mathbf{q_j}(\mathbf{x})$, собственные
значения которых равны $\mu_{j_0}$:
\begin{equation} \label{ker__1_}
Ker(\mathcal{N}_d-\mu_{j_0}\,I)= \sum_{\mu_j=\mu_{j_0}}
c_j\,\mathbf{q}_{j}(\mathbf{x}), \quad \text{для}\,\,\forall c_j \in
\mathcal{R}.\end{equation}
  \end{theorem}

%\newpage
\section {Построение собственных функций оператора $\nabla div$}
\subsection {Связь между собственными функциями операторов $\nabla div$
и Лапласа-Неймана в ограниченной области}
\begin{problem} Найти  все ненулевые собственные значения $\mu
$ и  собственные вектор-функции $\mathbf{u}(\mathbf{x})$ в
${{\mathbf{L}}_{2}}(G)$ оператора градиент дивергенции такие, что
\begin{equation}  \label{gd   1}-\nabla \text{ div }\mathbf{u}=
\mu \mathbf{u}\quad  \text {в} \quad G,\quad \mathbf{n}\cdot
\mathbf{u}{{|}_{\Gamma }}=0,
\end{equation}
              где
$\mathbf{n}\cdot \mathbf{u}$ - проекция вектора $\mathbf{u}$ на
нормальный вектор $\mathbf{n}$.\end{problem}
 К области определения
%$\mathcal{M}_{\mathcal{N}_d}$
 оператора $\mathcal{N}_d $ задачи 3 отнесем все вектор-функции
 $\mathbf{u}(\mathbf{x})$ класса $\mathcal{C}^2(G)\cap
 \mathcal{C}^1(\overline{G})$, которые удовлетворяют граничному условию
 $ \gamma\mathbf{n}\cdot\mathbf{u}=0$ и условию $\nabla \text{ div }\,\mathbf{u}\in
 {\mathbf{L}}_{2}(G)$.

Эта задача связана со спектральной задачей Неймана для скалярного
оператора Лапласа.

 \begin{problem} Найти все собственные значения
$\nu $ и собственные функции $g (\mathbf{x})$  оператора Лапласа
$-\Delta $ такие, что
      \begin{equation}  \label{gen   1}
              -\Delta g =\nu g \quad\text{в} \,\, G,\quad
 \mathbf{n}\cdot\nabla\,g|_{\Gamma }=0.
              \end{equation}\end{problem}
 К области определения
%$\mathcal{M}_{\mathcal{N}_{\Delta}}$
 оператора $\mathcal{N}_{\Delta}$ задачи 4 относят все функции
 $g(\mathbf{x})$ класса $\mathcal{C}^2(G)\cap
 \mathcal{C}^1(\overline{G})$, удовлетворяющие  условиям
 $\gamma \mathbf{n}\cdot\nabla
\,g=0$ , $\Delta\,g\in
 {{L}}_{2}(G)$. Эта задача является самосопряженной \cite{vla, mi}.
Решения задач 3 и 4 принадлежат классу
$\mathcal{C}^\infty(\overline{G})$, так как $\Gamma\in
\mathcal{C}^\infty$.

 Легко убедиться, что
 \begin{lemma} Любому решению
$(\mu ,\mathbf{u})$ задачи 3 в области G соответствует
 решение $(\nu ,g )=(\mu, \text{div }\mathbf{u})$ задачи 4.
  Обратно, любому решению  $(\nu ,g )$ задачи 4
соответствует решение  $(\mu ,\mathbf{u})=(\nu, \nabla g)$ задачи
3.\end{lemma}
    \subsection {Явные решения спектральной задачи Лапласа-Неймана в шаре }
   Согласно книге \cite{vla} В.С.Владимирова

    {\it собственные значения
    оператора Лапласа-Неймана $\mathcal{N}_\Delta$ в шаре  $B$   равны
    $\nu _{n,m}^{2}$,  где   $\nu
     _{n,m}^{{}}={{\alpha }_{n,m}}{{R}^{-1}}$,    $n\ge 0$,   $m\in N$, а   числа
      ${{\alpha }_{n,m}}>0$ суть нули  функций
      ${{\psi }_{n}^{\prime }}(z)$,
     производных ${{\psi }_{n}}(z)$, т.е. ${{\psi }_{n}}^{\prime }({{\alpha }_{n,m}})=0$.
Соответствующие $\nu _{n,m}^{2}$ собственные функции $g _{\kappa
}^{{}}$ имеют вид:
 \begin{equation}  \label{lan   1}g _{\kappa }^{{}}(r,\theta ,\varphi )=c{{_{\kappa }^{{}}}^{{}}}
 {{\psi }_{n}}{{(\alpha _{n,m}^{{}}r/R)}^{{}}}Y_{n}^{k}(\theta ,\varphi ),
  \end{equation}
где  $\kappa =(n, m, k)$- мультииндекс, \, $c_{\kappa
}^{{}}$-произвольные действительные постоянные,\,  $Y_{n}^{k}(\theta
,\varphi )$ - действительные сферические функции, \, $n\ge 0$,\,
$|k|\le n, \, m\in N$.}

Функции $g _{\kappa }^{{}}(x)$  принадлежат классу ${{C}^{\infty
}}(\overline{B})$  и при различных $\kappa$ ортогональны  в
${{L}_{2}}(B)$.   Система   функций $\{g _{\kappa }^{{}}\}$
полна в ${{L}_{2}}(B)$ \cite{mi}.
Нормируя их, получим  ортонормированный в ${{L}_{2}}(B)$ базис.

\subsection {Решение спектральной задачи 3 для $\nabla div$ в шаре}
Согласно лемме 3 вектор-функции ${{\mathbf{q}}_{\kappa }}(x)=
     \nabla {{g }_{\kappa }}(x)$
     являются решениями задачи 3 при ${\mu}_{n,m} ={\alpha }_{n,m}^2R^{-2}$ в
     ${{\mathbf{L}}_{2}}(B)$.
     Их компоненты $(q_r,q_\theta, q_\varphi)$ имеют вид
     \begin{equation}  \label{qom   1}\begin{array}{c}
     q _{r,\kappa }^{{}}(r,
     \theta ,\varphi )=c_{\kappa }(\alpha _{n,m}/R)
 {{\psi }_{n}^{{\prime}}}{{(\alpha _{n,m}^{{}}r/R)}^{{}}}Y_{n}^{k}
 (\theta ,\varphi ),\\
(q_{\varphi}+iq_{\theta})_{\kappa}=c_{\kappa }(1/r)
 {\psi }_{n}(\alpha _{n,m}r/R)\text{H}Y_{n}^{k}
 (\theta ,\varphi ).
  \end{array}\end{equation}
При $\kappa=(0,m,0)$ функция $Y_{0}^{0} (\theta ,\varphi )=1$,
$\text{H}Y_{0}^{0} =0$. Поэтому
\begin{equation}  \label{qomo   1}\begin{array}{l}
     q _{r,(0,m,0) }^{{}}(r)=c_{(0,m,0) }(\alpha _{0,m}/R)
 {{\psi }_{0}^{{\prime}}}{{(\alpha _{0,m}r/R)}},\\
(q_{\varphi}+iq_{\theta})_{(0,m,0)}=0.
  \end{array}\end{equation}

 Отметим, что $ {{\mathbf{q} }_{\kappa }}$ и $ {{\mathbf{q}
}_{{{\kappa }'}}}$ ортогональны при  ${\kappa }'\ne \kappa $.

Действительно, используя формулу Гаусса-Остроградского и свойства
этих векторов имеем

\begin{equation}  \label{ort   2}
 \int\limits_{B}{ {{\mathbf{q} }_{{{\kappa }'}}}\cdot  {{\mathbf{q}
}_{\kappa
}}^{{}}dx{{=}^{{}}}\frac{{\alpha}_{n,m}^2}{R^2}\int\limits_{B}{{{g
}_{{{\kappa }'}}}^{{}}{{g }_{\kappa }}^{{}}dx}}.\end{equation}
 Но функции
${{g}_{\kappa }}(x)$ и  ${{g }_{{{\kappa }'}}}(x)$
%, согласно(\ref{lan 1}),
взаимно ортогональны  в ${{\text{L}}_{2}}(B)$
 Значит,
   вектор - функции $ {{\mathbf{q} }_{\kappa }}$ и
$\mathbf{q}_{{\kappa }'}$ также взаимно ортогональны в
${{\mathbf{L}}_{2}}(B)$ и
 $\left\|
{{\mathbf{q}}_{\kappa }}(x) \right\|=({\alpha}_{n,m}/{R})\left\|
{{{g}}_{\kappa }}(x) \right\|$.

\subsection {Решение спектральной задачи  1 для ротора при $\lambda =0$ в шаре}
Числа $\mu_{n,m}={\alpha }_{n,m}^2
     R^{-2}>0$ при любых    $n\ge 0$,   $m\in N$.
     Поэтому вектор-функции $\mathbf{q}_{\kappa }$ являются также решениями
задачи 1 при $\lambda =0$.

 \subsection{Сходимость ряда Фурье по  собственным функциям
  оператора Лапласа-Неймана  в норме пространства Соболева}
В $\S\ 2.5$ главы 4 книги В.П.Михайлова \cite{mi} для областей $G$ с
границей $\Gamma\in \mathcal{C}^s$ определены подпространства
$H^s_\mathcal{N}(G)$  в $H^s(G)$:
\begin{equation} \label{sdp 1} H^s_\mathcal{N}(G)=\{f\in H^s(B):
(\mathbf{n}\cdot\nabla) f|_S=0,...,(\mathbf{n}\cdot\nabla)
  \triangle^{\sigma}f|_S=0\}, \end{equation}
где $\sigma$ равна целой части $[s/2-1]$ числа  $s/2-1$,
  $ s\geq 2$, и
    $H^0_\mathcal{N}(G)=L_2(G)$,
  $H^1_\mathcal{N}(G)=H^1(G)$ по определению.
Доказано, что принадлежность $f$ пространству $H^s_\mathcal{N}(G)$
необходима и достаточна для сходимости ее ряда Фурье по системе
собственных функций %${g}_{\kappa}$
оператора Лапласа-Неймана в $H^s(G)$(см. теоремы 8 и 9 $\S\ 2.5$ гл.
4).

%%%%%
\subsection{Сходимость ряда \eqref{arof 1}
  в норме пространства Соболева $H^s(B)$} Определим подпространство
   $\mathbf{A}^s_\mathcal{K}(B)$ в $\mathcal{A}$
при $s\geq 1$:
\[ \mathbf{A}^s_\mathcal{K}(B)= \{\mathbf{f}\in
\mathcal{A}\cap\mathbf{H}^s(B): \mathbf{n}\cdot\mathbf{f}|_S=0,...,
 \mathbf{n}\cdot\ (\nabla\text{div})^{\sigma}\mathbf{f}|_S=0,\,\,% \quad
 \|\mathbf{f}\|_{\mathbf{A}_\mathcal{K}^s}=
\|\mathbf{f}\|_{\mathbf{H}^s}\},\]
 где $\sigma=[(s-1)/2]+1$.
  Имеет место
\begin{theorem}  Для того, чтобы $\mathbf{f}\in
\mathcal{A}$
 разлагалась в ряд Фурье
 \begin{equation} \label{arof 1}
\mathbf{f}(\mathbf{x})=\sum_{\kappa}
(\mathbf{f},\mathbf{q}_{\kappa})\mathbf{q}_{\kappa}(\mathbf{x})
\end{equation}
 по системе собственных вектор-функций $\mathbf{q}_{\kappa}(\mathbf{x})$
 оператора градиента дивергенции в шаре,
 сходящийся в норме
 пространства Соболева $\mathbf{H}^s(B)$, необходимо и достаточно,
  чтобы $\mathbf{f}$ принадлежала $\mathbf{A}^s_\mathcal{K}(B)$.

 Если $\mathbf{f}\in \mathbf{A}^s_\mathcal{K}(B)$,
то сходится ряд
\begin{equation} \label{arof 2}
\sum_{\kappa}{\nu}_{\kappa}^{2s}\,
|(\mathbf{f},\mathbf{q}_{\kappa})|^2 ,\quad
{\nu}_{\kappa}=({\alpha}_{n,m})/R
\end{equation} и существует такая положительная постоянная $C>0$, не
зависящая от $\mathbf{f}$, что
\begin{equation} \label{aorf 3}
\sum_{\kappa} {\nu}_{\kappa}^{2s}\,
|(\mathbf{f},\mathbf{q}_{\kappa})|^2 \leq
C\|\mathbf{f}\|^2_{\mathbf{H}^s(B)}.
\end{equation}
  Если $s\geq 2$, то любая вектор-функция  $\mathbf{f}$ из
  $\mathbf{A}^s_\mathcal{K}(B)$
разлагается в в ряд Фурье, сходящийся в пространстве $\mathbf{C}^{s-2}(\overline{B})$.%
\end{theorem}
Доказательство этой теоремы аналогично  доказательству
 Теоремы 3 для оператора ротор.

{\bf Следствие.} {\it  Любая вектор-функция $f$ из
$\mathcal{A}\cap\mathbf{C}^{\infty}_0({B})$ разлагается в
 ряд Фурье \eqref{arof 1},
сходящийся в пространстве $\mathbf{C}^{\infty}(\overline{B})$.}

\subsection{Скалярное произведение функций $\mathbf{f}$ и
$\mathbf{g}$ из $\mathcal{A}_{\gamma}$ в базисе из собственных
функций градиента дивергенции} Оно имеет вид:
\begin{equation} \label{spa 2}( \mathbf{f}, \mathbf{g})=
\sum_{\kappa, n\geq 0} \,
(\mathbf{f},\mathbf{q}_{\kappa})(\mathbf{g},\mathbf{q}_{\kappa})
\end{equation}
Если $\mathbf{f}$ и $\mathbf{g}$ принадлежат
$\mathbf{A}^1_\mathcal{K}(B)$, то равенства
\begin{equation} \label{sgd }(\nabla\bf{div}\, \mathbf{f}, \mathbf{g})=
(\mathbf{f},\nabla\text{div}\mathbf{g})= \sum_{\kappa, n\geq 0}
{\nu}_{\kappa}^2[(\mathbf{f},\mathbf{q}_{\kappa})
(\mathbf{g},\mathbf{q}_{\kappa}) ]
\end{equation}показывают, что оператор $\nabla\mathbf{div}$ является
самосопряженным в $\mathcal{A}_{\gamma}$.

\section{Базисные подпространства в $\mathbf{L}_{2}(G)$}

\subsection{Пространства $\mathcal{A}$ и $\mathcal{B}$ в
${{\mathbf{L}}_{2}}({G})$} В $\S 1$ мы рассмотрели ортогональные
подпространства $\mathcal{A}$, $\mathcal{B}$, $\mathcal{B}_H$ и
$\mathbf{V}^0$. Пространства $\mathcal{A}$ и $\mathcal{B}_H$
принадлежат ядру оператора ротор; а в $\mathbf{V}^0$ он продолжается
как самосопряженный оператор $S$,  собственные функции которого
образуют %ортонормированный
 базис в $\mathbf{V}^0$.% Согласно п. 4.3

Пространство $\mathcal{B}$ %=\mathcal{B}_H\oplus\ \mathbf{V}^0$
принадлежат ядру оператора градиент дивергенции; а в
$\mathcal{A}_\gamma$ он продолжается как самосопряженный оператор
$\mathcal{N}_d$, собственные функции которого образуют
%ортонормированный
базис в $\mathcal{A}_\gamma\subset \mathcal{A}$.

 $\mathcal{B}_H$--конечномерное пространство.

Согласно \eqref{wor  1} и \eqref{bhr  1} пространство ${\mathbf{{L}}_{2}}(G)$ разлагается на ортогональные
подпространства: % собственным вектор-функциям оператора ротор.
\begin{equation}
\label{wro 1}\mathbf{L}_{2}(G)=\mathcal{A}\oplus\mathcal{B}_H\oplus
\mathbf{V}^{0}(G).
\end{equation}
Следовательно, имеет место
\begin{theorem} Система
$\{\mathbf{q}_{j}(x)\}\cup\{\mathbf{h}_j(x)\}\cup\{\mathbf{q}_{j
}^{+}(x)\}\cup \{\mathbf{q}_{j }^{-}(x)\}$ решений задачи 1 образует
в пространстве ${\mathbf{{L}}_{2}}(G)$ ортонормированный базис.
Любую вектор-функцию  $\mathbf{f}(\mathbf{x})$ из
${\mathbf{{L}}_{2}}(G)$
 можно разложить в ряд Фурье по этому базису:
\[\mathbf{f}(x)=\sum_{\mu_j\in M}(\mathbf{f},\mathbf{q}_{j})\mathbf{q}_{j}
(x)+ \sum_{j=1}^{\rho}(\mathbf{f},\mathbf{h}_{j})
 \mathbf{h}_{j}(x)+\sum_{\lambda_j\in
\Lambda}[(\mathbf{f},\mathbf{q}^{+}_{j})\mathbf{q}^{+}_{j}+
(\mathbf{f},\mathbf{q}^{-}_{j})\mathbf{q}^{-}_{j}].
\]
 \end{theorem}

В случае шара пространство $\mathcal{B}_H$ пусто.
\subsection{Разложение векторного поля $\mathbf{f}\in \mathbf{L}_{2}(B)$ }
% Это разложение   из
 на безвихревое поле  $\mathbf{a}_f$
 и соленоидальное поле  $\mathbf{b}_f$:\quad
$\mathbf{f}(\mathbf{x})=\mathbf{a}_f(\mathbf{x})+\mathbf{b}_f(\mathbf{x})$,
где
\begin{equation} \label{sra 1}
{\mathbf{a}_f}=%\mathcal{P}_{\mathcal{A}}f\equiv
\sum_{n=0}^{\infty}\sum_{m=1}^{\infty}\sum_{k=-n}^{n}
(\mathbf{f},{\mathbf{q}}_{\kappa})\,\mathbf{q}_{\kappa}(\mathbf{x}),\quad
\kappa=(n,m,k)
\end{equation}

\begin{equation} \label{srpm 2}
{\mathbf{b}_f}=\sum_{n=1}^{\infty}\sum_{m=1}^{\infty}\sum_{k=-n}^{n}
[(\mathbf{f},{\mathbf{q}}_{\kappa}^+)\,\mathbf{q}_{\kappa}^+(\mathbf{x})
+(\mathbf{f},{\mathbf{q}}_{\kappa}^-)\,\mathbf{q}_{\kappa}^-(\mathbf{x})].
\end{equation}%\end{array}
Частичные суммы $\mathbf{S}^0_N$ и $\mathbf{S}^1_N$ рядов \eqref{sra
1} и \eqref{srpm 2} состоят из элементов с индексами $n,m,k$, для
которых $0<\alpha_{n,m}<N$ и
 $0<\rho_{n,m}<N$, соответственно.% а затем $N\rightarrow\infty$.

 Имеет место
равенство Парсеваля-Стеклова:
$\|\mathbf{f}\|^2=\|\mathbf{a}_f\|^2+\|\mathbf{b}_f\|^2$, которое
запишем так
\begin{equation} \label{psra 1}
\|{\mathbf{f}}\|^2=\sum_{N=1}^{\infty} \sum_{(n,m)\in
\mathbb{P}_N}^{}\sum_{k\in[-n,n]}^{}
[(\mathbf{f},{\mathbf{q}}_{\kappa})^2
+(\mathbf{f},{\mathbf{q}}_{\kappa}^+)^2
+(\mathbf{f},{\mathbf{q}}_{\kappa}^-)^2],
\end{equation}где решетка
 $\mathbb{P}_N=\{(n,m):0<\rho _{n,m}<N,\, 0<\alpha_{n,m}<N\}$ и векторы
${\mathbf{q}}_{0,m,0}^{\pm}=0$. % (см. $\S$ 2).

Отметим, что разложение векторного поля $\mathbf{f}(\mathbf{x})$ на
безвихревое поле $\nabla{h}(\mathbf{x})$ и соленоидальное поле
$\mathbf{u}(\mathbf{x})$ связано  с решением задачи Неймана
\begin{equation} \label{npr 1}
\triangle h= \text{div}\,\mathbf{f} \quad \text{в}\quad B,\quad
\mathbf{n}\cdot\nabla\,h|_S=\mathbf{n}\cdot\mathbf{f}|_S,
\end{equation}в классической или обобщенной постановках
(cм.\cite{lad, bs}).

 Мы  сводим
 решение  задачи к вычислению интегралов
$(\mathbf{f},{\mathbf{q}}_{\kappa})$,
$(\mathbf{f},{\mathbf{q}}_{\kappa}^+)$,
$(\mathbf{f},{\mathbf{q}}_{\kappa}^-)$
 и ее решение получаем в виде рядов
\eqref{sra 1}, \eqref{srpm 2}.

\section{Решение  краевой задачи  5 в шаре}
 Методом Фурье легко
решается  краевая
\begin{problem} Пусть задана вектор-функция $\mathbf{f}(\mathbf{x})\in
\mathbf{L}_{2}(B)$. Найти вектор-функцию $\mathbf{u}(\mathbf{x})$ в
$\mathbf{H}^{1}(B)$  такую, что
 \begin{equation}
\label{rot__2_} {\bf rot}\, \mathbf{u}+\lambda \mathbf{u}=
\mathbf{f}\quad \text{в}\quad B,\quad \gamma\mathbf{n}\cdot
\mathbf{u}=0.\end{equation}
\end{problem}

\subsection {Основные пространства}
Через $ \mathbf{E}^{s}(B)$ %(или $\mathbf{H}^{s}_{div}(B)$)
обозначают \cite{rt}
пространство% в $\mathbf{L}_{2}(B)$:
\[\mathbf{E}^{s}(B)=\{\mathbf{v}\in
\mathbf{H}^{s}(B):\,\text{div}\mathbf{v}\in
{H}^{s}(B),\,\,\|\mathbf{v}\|_{\mathbf{E}^{s}}=(
 \|\mathbf{v}\|^2_{\mathbf{H}^{s} }+
 \|\text{div}\mathbf{v}\|^2_{{H}^{s}})^{1/2}\},\]
  где число
 $s\geq 0$ целое.   Оно является
 пространством Гильберта и
  \begin{equation} \label{vlo s}
\mathbf{C}_0^{\infty}({B})  \subset \mathbf{E}^{s}(B), \quad
 \mathbf{H}^{s+1}(B)\subset \mathbf{E}^{s}(B)
 \subset\mathbf{H}^s(B).  \end{equation}
Согласно п. 2.2, ${\bf rot}\, \mathbf{u}+\lambda \mathbf{u}\in
\mathbf{E}^{s}(B)$, если $\mathbf{u}\in \mathbf{H}^{s+1}(B)$.

Как известно\,\cite{mi},  для функций $v$ из пространства $H^1(B)$
определен оператор {\it следа} $\gamma :H^1(B)\rightarrow
H^{1/2}(S)$, равный следу $v$ на $S$ для гладких функций из
$\mathcal{C}^1(\overline{B})$: $\gamma\,v=v|_S$,
 причем $\quad \|\gamma\,v\|_{L_2(S)}\leq c \|v\|_{H^1(B)}$.

 Аналогично, для вектор-функций $\mathbf{u}(\mathbf{x})$ из
$\mathbf{E}^0(B)$ определен \cite{rt} оператор {\it следа нормальной
компоненты} $\gamma_\mathbf{n} :\mathbf{E}^0(B)\rightarrow
H^{-1/2}(S)$, равный сужению $\mathbf{n}\cdot \mathbf{u}$ на $S$ для
функций из
$\mathcal{C}^1(\overline{B})$:
$\gamma_{\mathbf{n}}\,\mathbf{u}={\mathbf{n}}\cdot\mathbf{u}|_S$.

 Для $u\in \mathbf{E}^0(B)$ и $v\in H^1(B)$ верна
обобщенная формула Стокса:
%\begin{equation} \label{ofc 1}
 $\langle \gamma_{\mathbf{n}}\,{\mathbf{u}},\gamma\,v \rangle=
 (\mathbf{u}, \nabla v)+( \text{div}\mathbf{u}, v)$
 %\end{equation}
 где  $\langle\gamma_{\mathbf{n}}\,{\mathbf{u}},\gamma\,v \rangle$- линейный
функционал над %пространством
$H^{1/2}(S)$; $\gamma_{\mathbf{n}}\,{\mathbf{u}}\in H^{-1/2}(S)$ при
$\gamma\,v\in H^{1/2}(S) $. Имеют место непрерывные вложения:
 %\begin{equation} \label{vlo 1}
 $ H^{1/2}(S)\subset L_2(S)\subset H^{-1/2}(S)$
%\end{equation}

Пространство
$\mathbf{{E}}_{\gamma}^s(B)=
 \{\mathbf{f}\in\mathbf{E}^s(B):
{\mathbf{n}}\cdot{\mathbf{f}}|_S=0 \}, \quad s\geq 0.$

 \subsection {Решение  краевой задачи  \eqref{rot__2_}
при  $\lambda\neq Sp\, (rot)$}
 \begin{theorem} Если $\lambda\neq 0, \pm {\lambda}_{n,m}$,\,
 $n,m\in \mathbf{N}$ и
  $\mathbf{f}\in\mathbf{E}_{\gamma}^{0}(B)$,
    то
 единственное решение $\mathbf{u}$ задачи 5 дается суммой  рядов
 $\mathbf{u}_1+\mathbf{u}_2$, где
 \begin{equation} \label{kra 1}
{\mathbf{u}_1}={\lambda}^{-1}
\sum_{n=0}^{\infty}\sum_{m=1}^{\infty}\sum_{k=-n}^{n}
(\mathbf{f},{\mathbf{q}}_{\kappa})\,\mathbf{q}_{\kappa}
(\mathbf{x}),\quad \kappa=(n,m,k),
\end{equation}
\begin{equation} \label{kro 2}{\mathbf{u}_2}=
\sum_{n=1}^{\infty}\sum_{m=1}^{\infty}\sum_{k=-n}^{n}
[\frac{(\mathbf{f},{\mathbf{q}}_{\kappa}^+)}{\lambda+\lambda_{n,m}}\,
\mathbf{q}_{\kappa}^+(\mathbf{x})+
\frac{(\mathbf{f},{\mathbf{q}}_{\kappa}^-)}{\lambda-\lambda_{n,m}}\,
\mathbf{q}_{\kappa}^-(\mathbf{x})].
\end{equation}
Решение задачи принадлежит пространству Соболева
$\mathbf{H}^{1}_\gamma(B)$.\newline Если $\mathbf{f}\in
\mathcal{A}\subset \mathbf{L}_{2}(B)$, то
  $\mathbf{u}={\lambda}^{-1}\mathbf{f}$ отображает $\mathcal{A}$
  на $\mathcal{A}$.\newline
Если    $\mathbf{f}\in \mathcal{B}\bot \mathcal{A}$ в
 $\mathbf{L}_{2}(B)$, то
 $\mathbf{u}=\mathbf{u}_2$ принадлежит % отображает $\mathcal{B}$ на
$\mathbf{W}^{1}(B)\subset \mathbf{H}^{1}_{\gamma}(B)$.\newline Если
же $\mathbf{f}\in \mathcal{D}(B)$, то  ряды \eqref{kra 1},
\eqref{kro 2} сходятся в любом из пространств $\mathbf{H}^{s}(B)$,\,
$s\geq 1$ и их сумма есть  классическое решение задачи класса
$C^{\infty}(\overline{B})$.
\end{theorem}
Доказательство приведено в \cite{saUMJ13}.

 Мы не будем выписывать подробно  решение задачи 6 при $\lambda=0$.
 Отметим только, что условие $\text{div}\,\mathbf{f}=0$ необходимо
 и достаточно для ее разрешимости, а однородная задача имеет счетное
 число линейно независимых решений $\mathbf{q}_{\kappa}(\mathbf{x})$.

  При $\lambda=\pm\lambda_{n,m}$ задача \eqref{rot__2_}
разрешима по Фредгольму.

  Задача:
$\nabla\mathbf{div}\mathbf{v}+\lambda\mathbf{v}=\mathbf{f}$ в $B$,
$\gamma\mathbf{n}\cdot\mathbf{v}=0$ решается аналогично.
%%%%%%%

\newpage
Реферат:

Ряды Фурье оператора ротор и пространства Соболева.

  Р.С.Сакс

Автор
 изучает  свойства
 операторов ротор и  градиент дивергенции  в произвольной ограниченной
  области $G$
 с гладкой границей  $\Gamma$,
  их спектральные разложения и краевые задачи.

% в
В пространстве $\mathbf{L}_{2}(G)$ выделяются ортогональные
подпространства $\mathbf{V}^{0}(G)$   и
${\mathcal{{A}}_{\gamma}}(G)$. Доказано, что существуют продолжения
$S$ и $\mathcal{N}_d$ \,  операторов ротор и градиент дивергенции в
эти пространства такие, что $S$ и $\mathcal{N}_d$ являются
самосопряженными и имеют вполне непрерывные обратные $S^{-1}$ и
$\mathcal{N}_{d}^{-1}$ . Откуда вытекает ортогональность собственных
функций каждого из этих операторов в соответствующих
подпространствах и полнота совокупной системы в  их объединении.
%${\mathbf{{L}}_{2}}(G)$.

Найдены необходимые и достаточные условия на $\mathbf{u}\in
\mathbf{V}^{0}(B)$ и $\mathbf{v}\in {\mathcal{{A}}_{\gamma}}(B)$,
при которых их ряды Фурье сходятся в норме пространства Соболева
$\mathbf{H}^{s}(B)$.

При  $\lambda\neq 0$ исследована разрешимость в $\mathbf{H}^{s}(G)$
краевых задач для систем
$\mathbf{rot}\mathbf{u}+\lambda\mathbf{u}=\mathbf{f}$ в  $G$ и
$\nabla\mathbf{div}\mathbf{u}+\lambda\mathbf{u}=\mathbf{f}$ в  $G$ с
условием $\mathbf{n}\cdot\mathbf{u}|_{\Gamma}=g$ на границе.

Методом Фурье при любых $\lambda$ решена краевая задача для системы
$\mathbf{rot}\mathbf{u}+\lambda\mathbf{u}=\mathbf{f}$ в шаре $B$ с
условием $\mathbf{n}\cdot\mathbf{u}|_S=0$.

 Fourier series of the curl  operator and Sobolev  spaces

Saks Romen Semenovich

Сакс Ромэн Семенович ведущий научный сотрудник
 Институт Математики с
ВЦ УНЦ РАН 450077, г. Уфа, ул. Чернышевского, д.112 телефон: (347)
272-59-36
                 (347) 273-34-12
факс:        (347) 272-59-36 телефон дом.: (347) 273-84-69 моб.
+79173797538

 e-mail: romen-saks@yandex.ru

 \end{document}